\documentclass[11pt]{article}

\usepackage{lscape}
\usepackage{epsfig}
\usepackage{enumerate}
\usepackage{amssymb}
\usepackage{afterpage}
\usepackage{theorem}
\usepackage{xypic}

\newcommand{\R}{\mbox{$\mathbb{R}$}}

\newcommand{\C}{\mbox{$\mathbb{C}$}}

\newcommand{\D}{\mbox{\bf D}}

\newcommand{\mod}{\mbox{mod}\;}

\newcommand{\ov}{\overline}

\newcommand{\dfrac}{\displaystyle\frac}
\newcommand{\proof}{\noindent{\bf Proof:} \quad}
\newcommand{\proofof}[1]{\noindent {\bf Proof of #1} \hspace{0.1in}}
\newcommand{\qed}{\hfill\mbox{\raggedright\rule{0.07in}{0.1in}}\vspace{0.1in}}

\newcommand{\dsum}{\displaystyle\sum}

\newtheorem{theorem}{Theorem}[section]
\newtheorem{proposition}[theorem]{Proposition}
\newtheorem{corollary}[theorem]{Corollary}
\newtheorem{lemma}[theorem]{Lemma}

\theorembodyfont{\rmfamily}

\newtheorem{example}[theorem]{Example}
\newtheorem{remark}[theorem]{Remark}

\newcommand{\diag}{\mbox{{\rm diag}}}

\newcommand{\AND}{\quad\mbox{and}\quad}

\title{Realization of critical eigenvalues for scalar and symmetric linear delay-differential equations}
\author{P-L. Buono\\
Faculty of Science\\
University of Ontario Institute of Technology\\
Oshawa, ONT L1H 7K4\\
Canada\\
\\
V.G. LeBlanc\\
Department of Mathematics and Statistics\\
University of Ottawa\\
Ottawa, ONT K1N 6N5\\
Canada}
\date{\today}
\begin{document}
\maketitle

\begin{abstract}
This paper studies the link between the number of critical
eigenvalues and the number of delays in certain classes of
delay-differential equations. There are two main results. The first
states that for $k$ purely imaginary numbers which are linearly
independent over the rationals, there exists a scalar
delay-differential equation depending on $k$ fixed delays whose
spectrum contains those $k$ purely imaginary numbers. The second
result is a generalization of the first result for
delay-differential equations which admit a characteristic equation
consisting of a product of $s$ factors of scalar type. In the second
result, the $k$ eigenvalues can be distributed amongst the different
factors. Since the characteristic equation of scalar equations
contain only exponential terms, the proof exploits a toroidal
structure which comes from the arguments of the exponential terms in
the characteristic equation. Our second result is applied to delay
coupled $\D_n$-symmetric cell systems with one-dimensional cells. In
particular, we provide a general characterization of delay coupled
$\D_n$-symmetric systems with arbitrary number of delays and cell
dimension.
\end{abstract}

\section{Introduction and Background}

Delay-differential equations (DDEs) have been used as mathematical
models for phenomena in population dynamics~\cite{Kuang},
physiology~\cite{GM,BGMT-book}, physics~\cite{LK}, climate
modelling~\cite{SS} and engineering~\cite{SS} amongst others.
Delay-differential equations behave like abstract ordinary
differential equations (ODEs) on an infinite-dimensional (Banach)
phase space and many results which are known for ODEs on
finite-dimensional spaces have analogs in the context of DDEs. Many
scalar delay-differential equation models have been developed over
the years such as for Cheynes-Stokes respiration~\cite{GM} and the
regulation of hematopoiesis~\cite{GM}, the delayed Nicholson
blowflies equation~\cite{Gurney} in population dynamics, a two-delay
model of an experiment on Parkinsonian tremor~\cite{BBL} and many
more.

The bifurcation analysis of DDEs is done essentially in the same way
as in ODEs, although the technical details differ. Consider the
neighborhood of an equilibrium solution of a nonlinear DDE, then the
analysis of the linearization at the equilibrium point leads to
stable, unstable and centre invariant subspaces where only the
stable subspace is infinite-dimensional. There exists local
invariant manifolds (stable, unstable and center manifolds) tangent
to the corresponding invariant subspaces of the linearized equations
about the equilibrium point on which the flow near the equilibrium
is either exponentially attracting (stable manifold), exponentially
repelling (unstable manifold), or non-hyperbolic (center manifold).
Now, bifurcations near equilibria are determined by the flow on the
centre manifold and the dimension of this manifold is determined by
the number of eigenvalues of the linearization on the imaginary
axis.

The first result of our paper is Theorem~\ref{thm:main} and goes as
follows. Consider $n$ nonzero imaginary numbers
$i\omega_1,\ldots,i\omega_n$ where the imaginary parts
$\omega_1,\ldots,\omega_n$ are positive and not rationally
dependent. We show that there exists a scalar linear
delay-differential equation depending on $n$ discrete delays written
\begin{equation}\label{linear-eq1}
\dot x=\sum_{j=1}^{n} a_j x(t-\tau_j)
\end{equation}
where $x\in \R$, $a_j\in \R$ and $\tau_j\in [0,\tau]$ for all
$j=1,\ldots,n$ such that the characteristic equation
of~(\ref{linear-eq1}),given by
\begin{equation}\label{char-eq}
\lambda-\sum_{j=1}^{n} a_j e^{\lambda\tau_j}=0,
\end{equation}
has eigenvalues $\pm i\omega_1,\ldots,\pm i\omega_n$. This result
generalizes explicit computations done in the case of one and two
delays, see~\cite{HL,Diekmann-etal,BC94}. The proof is done by
embedding the problem as a mapping which is solved by the implicit
function theorem at a carefully chosen point. From the implicit
function theorem, we are able to define a smooth mapping whose
transversal intersection with a dense curve on an $n$-dimensional
torus provide solutions. The incommensurability of the $n$
frequencies enables us to define the dense curve on the $n$-torus.
This type of argument using a dense curve on an $n$-dimensional
torus was used in Choi and LeBlanc~\cite{Choi-LeBlanc1}.

This result falls within the category of so-called realization
theorems. For instance, the realization theorem of linear ODEs by
linear DDEs obtained by Faria and Magalh$\tilde{\mbox{\rm
a}}$es~\cite{FMR}. They show that for any finite dimensional matrix
$B$, a necessary and sufficient condition for the existence of a
bounded linear operator ${\cal L}_0$ from $C([-\tau,0],{\mathbb
R}^{n})$ into ${\mathbb R}^n$ with infinitesimal generator having
spectrum containing the spectrum of $B$ is that $n$ be larger than
or equal to the largest number of Jordan blocks associated with each
eigenvalue of $B$. Other results in this direction are concerned
with the realization of finite jets of ODEs on a finite-dimensional
centre manifold by delay-differential equations,
see~\cite{FMR,Choi-LeBlanc1}. To our knowledge, the realization
theorems in this paper are the first general results linking the
number of critical eigenvalues of linear delay-differential
equations with the number of discrete delays.

The next significant result is an openness theorem, that is, the
realization of $n$ imaginary numbers (not necessarily rationally
independent) as eigenvalues of a linear scalar delay-differential
equation is valid in a neighborhood of any set of $n$ rationally
independent imaginary numbers. The proof of this theorem also relies
on the implicit function theorem.

We then turn our attention to the context of symmetric systems of
delay-differential equations. Several examples of symmetric systems
of DDEs~\cite{Guo-Huang1,Peng07} have characteristic equations which
decompose in factors, some of which have the same form as the
characteristic equation~(\ref{char-eq}). The decomposition of the
characteristic equation is induced by the isotypic decomposition of
the space and we present a general derivation of this decomposition.
We show that isotypic components consisting of a unique
one-dimensional complex irreducible representation contribute a
factor of the form~(\ref{char-eq}) in the characteristic equation
and so Theorem~\ref{thm:main} can be applied directly to each such
factors separately.

We present a generalization of Theorem~\ref{thm:main} to the case
where several factors of the characteristic equation have purely
imaginary eigenvalues simultaneously. Theorem~\ref{thm:main2} shows
that a set of $n$ rationally independent purely imaginary complex
numbers can be realized from several factors of the characteristic
equation of a delay-differential equation with $n$ delays given some
nondegeneracy conditions on the characteristic equation are
satisfied. The statement of the theorem is independent of any
symmetric structure and the proof is a generalization of the proof
of Theorem~\ref{thm:main}.

We illustrate the above result on $\D_n$-symmetric rings of $n$
delay equations with delayed coupling. Hopf bifurcation from such
symmetric networks have been studied by several
authors~\cite{Campbelletal05,Guo05,Guo-Huang1,Guo-Huang2,Peng07,Peng-Yuan1,Wu98}.
In order to apply Theorem~\ref{thm:main2} to this context, we derive
an explicit form of the coupling matrix in terms of the connections
in the graph representation of the ring for cells of any dimension
and arbitrary numbers of connections and delays. This is a
generalization of the networks considered in the articles listed
above in this paragraph.

We specialize to the case of one-dimensional cells and we can then
obtain general formulae for the factors of the form~(\ref{char-eq})
in the characteristic equation of the $\D_n$-symmetric ring. Given a
nondegeneracy assumption of Theorem~\ref{thm:main2} is satisfied,
then it can be applied directly to $\D_n$ symmetric coupled cell
systems with $n$ odd. We do not treat the case $n$ even because the
form of the equations needed for Theorem~\ref{thm:main2} is not
satisfied. We illustrate this fact in a $\D_4$-symmetric example.


The paper is organized as follows. The first section contains brief
preliminary remarks and then we state and prove our main result
(Theorem~\ref{thm:main}) and the openness result
(Theorem~\ref{thm:open}). Then we introduce the context leading to
Theorem~\ref{thm:main2} and state this result.
Section~\ref{sec:symmetry} is devoted to $\Gamma$-symmetric systems
of delay-differential equations and the section begins with a
general discussion. Section~\ref{section:delay-ring} presents a
characterization of $\D_n$-symmetric rings of delay coupled cells
with an arbitrary number of delays and derive the characteristic
equation in the case of one-dimensional cells. In
Section~\ref{section:critical-Dn}, Theorem~\ref{thm:main2} is
applied to $\D_n$-symmetric rings of one-dimensional cells with $n$
odd. Section~\ref{sec:main-proof} has the proof of
Theorem~\ref{thm:main2}. We conclude by a discussion of open
problems along the lines of the ones presented in this paper.

\section{Realization theorems}

We now discuss some aspects of the spectral theory of linear scalar
delay-differential equations. In fact, we just introduce the basic
facts, in a non-abstract setting, needed for the statement of our
first main theorem. For a complete treatment, see Diekmann {\em et
al}~\cite{Diekmann-etal} or Hale et Lunel~\cite{HL}.

Consider the scalar delay-differential equation
\begin{equation}\label{eq:scalar}
\dot x(t)=\sum_{j=1}^{n} a_j x(t-\tau_j)
\end{equation}
where $a_j\in \R$ and $\tau_j\in [0,\tau]$ for all $j=1,\ldots,n$
and $\tau>0$. The characteristic equation for~(\ref{eq:scalar}) can
be obtained by substituting $x(t)=Ce^{\lambda t}$, where $C$ is a
constant, into the equation. Thus,
\[
\begin{array}{rcl}
\lambda Ce^{\lambda t}&=&\sum_{j=1}^{n} a_j
Ce^{\lambda(t-\tau_j)}=\sum_{j=1}^{n} a_j
Ce^{-\lambda\tau_j}e^{\lambda t}
\end{array}
\]
and by rearranging the terms we obtain
\[
\left(\lambda-\sum_{j=1}^{n} a_j e^{-\lambda\tau_j}\right)x(t)=0.
\]
So, $x(t)$ is a nonzero solution of~(\ref{eq:scalar}) if and only if
\[
\Delta(\lambda):=\lambda-\sum_{j=1}^{n} a_j e^{-\lambda\tau_j}=0.
\]
The complex number $\lambda$ is an eigenvalue of
equation~(\ref{eq:scalar}) if it is a solution of the characteristic
equation $\Delta(\lambda)=0$.

The question we address in this paper is related to the number of
imaginary eigenvalues (with incommensurable frequencies) which can
satisfy $\Delta(\lambda)=0$. The case $n=1$ with one nonzero delay
is a straightforward calculation and $\Delta(\lambda)=0$ for only
one nonzero imaginary eigenvalue $\lambda$, see~\cite{HL}. The case
$n=2$ with $\tau_1=0$ and $\tau_2\in (0,\tau]$ in~(\ref{eq:scalar})
can be found in~\cite{Diekmann-etal}. There, it is shown that
$\Delta(\lambda)=0$ can have at most two nonzero imaginary
eigenvalues. The case $n=2$ with $\tau_1,\tau_2> 0$ is done
in~\cite{BC94} where it is shown that $\Delta(\lambda)=0$ can have
at most two nonzero imaginary eigenvalues. We are now ready to state
our first result.
\begin{theorem}\label{thm:main} Suppose
$\omega_1>0,\,\omega_2>0,\,\ldots\,,\omega_{n}>0$ are linearly
independent over the rationals.  Then there exists
$\tau_1>0,\,\tau_2>0,\,\ldots\,,\tau_n>0$,
$a_1\in\mathbb{R},\,a_2\in\mathbb{R},\,\ldots\,,a_n\in\mathbb{R}$
such that the linear delay differential equation
\begin{equation}
\dot{x}(t)=a_1\,x(t-\tau_1)+a_2\,x(t-\tau_2)+\cdots+a_n\,x(t-\tau_n)
\label{DDE1}
\end{equation}
has solutions $x^{\pm}_j(t)=e^{\pm\,i\omega_jt}$ for all
$j=1,\ldots,n$.
\end{theorem}

\proof A necessary and sufficient condition for the conclusion of
the theorem to hold is that the following algebraic system of $2n$
equations
\begin{equation}
\begin{array}{ccc}
{\displaystyle\sum_{k=1}^n\,a_k\,e^{-i\omega_j\tau_k}}&=&i\omega_j,\,\,\,\,\,\,j=1,\ldots,n\\[0.15in]
{\displaystyle\sum_{k=1}^n\,a_k\,e^{i\omega_j\tau_k}}&=&-i\omega_j,\,\,\,\,\,\,j=1,\ldots,n
\end{array}
\label{chareq}
\end{equation}
has a solution in the $2n$ unknowns
$(\tau_1,\tau_2,\ldots,\tau_n,a_1,a_2,\ldots,a_n)$. Although
(\ref{chareq}) is in complex form, since the second equation in
(\ref{chareq}) is just the complex conjugate of the first equation
in (\ref{chareq}), system (\ref{chareq}) is equivalent to a system
of $2n$ {\it real} equations. This fact is taken for granted
throughout the sequel, even though we continue to use complex
notation.

It is useful to use the following matrix notation for (\ref{chareq})
\begin{equation}
\left(\begin{array}{c}P(\tau;\omega)\\[0.1in]P(-\tau;\omega)\end{array}\right)A^T=\left(\,\begin{array}{c}i\omega^T\\[0.1in]-i\omega^T\end{array}\,\right)
\label{chareq_matrix}
\end{equation}
where, $\omega=(\omega_1,\ldots,\omega_n)$, $A=(a_1,\ldots,a_n)$,
superscript $T$ denotes transpose, and
$P(\tau;\omega)=P(\tau_1,\ldots,\tau_n;\omega_1,\ldots,\omega_n)$ is
the $n\times n$ matrix whose entry at row $j$ column $k$ is
\[
\left[\,P(\tau;\omega)\,\right]_{jk}=e^{-i\omega_j\tau_k}.
\]
Note that $\overline{P(\tau;\omega)}=P(-\tau;\omega)$.

\vspace*{0.25in} \noindent {\bf Remark}: {\it Let ${\cal P}$ denote
the $2n\times n$ matrix of coefficients of the left-hand side of
(\ref{chareq_matrix}). If we view
$\tau_1\in\mathbb{R}^+,\ldots,\tau_n\in\mathbb{R}^+$ as free
parameters for these coefficients, we see that each column of the
matrix ${\cal P}$ is dense on an $n$-torus, $\mathrm{T}\equiv (\mathbb{S}^1)^n$, since
the $\omega_j$ are rationally independent. This remark is a key
point in the following discussion.}

\vspace*{0.25in}
Define $\mathrm{V}\equiv\mathrm{T}^n=((\mathbb{S}^1)^n)^n$, and choose
coordinates on $\mathrm{V}$ as follows:
\[
\mathrm{V}=\{\,\Phi=(\Phi^1,\ldots,\Phi^n)\,|\,\Phi^j=(\varphi^j_1,\ldots,\varphi^j_n)\in
(\mathbb{S}^1)^n,\,j=1,\ldots,n\,\}.
\]
Consider the following mapping associated to
(\ref{chareq_matrix}):
\[
F:\mathrm{V}\times\mathbb{R}^n\longmapsto\mathbb{R}^{2n}
\]
defined by
\begin{equation}
F(\Phi,A;\omega)=\left(\begin{array}{c}\widetilde{P}(\Phi)\\[0.1in]
\widetilde{P}(-\Phi)\end{array}\right)A^T-i\left(\begin{array}{c}\omega^T\\[0.1in]-\omega^T\end{array}\right),
\label{real_suspended}
\end{equation}
where $A$ and $\omega$ are as previously defined,
and $\widetilde{P}(\Phi)$ is the $n\times n$ matrix
whose entry at row $j$ column $k$ is
\[
\left[\,\widetilde{P}(\Phi)\,\right]_{jk}=\left[\,\widetilde{P}(\Phi^1,\ldots,\Phi^n)\,\right]_{jk}=e^{-i\varphi^k_j}.
\]

For reasons that are clear from the sequel, we choose to single out
the last torus in the product $\mathrm{V}=\mathrm{T}^n$, and
use the variable $\Psi$ as coordinates for this torus, i.e.
\[
\Psi=(\psi_1,\ldots,\psi_n)=\Phi^n=(\varphi^n_1,\ldots,\varphi^n_n).
\]
Thus, we write $\mathrm{V}=\mathrm{V}_{\Phi}\times\mathrm{V}_{\Psi}$,
where $\mathrm{V}_{\Phi}\cong \mathrm{T}^{n-1}$ and
$\mathrm{V}_{\Psi}\cong\mathrm{T}$, so that
\[
F:\mathrm{V}_{\Phi}\times\mathrm{V}_{\Psi}\times\mathbb{R}^n\longmapsto
\mathbb{R}^{2n}
\]
is written as $F(\Phi,\Psi,A;\omega)$ in (\ref{real_suspended}) (we
have re-labeled $\Phi=(\Phi^1,\ldots,\Phi^{n-1})$ to designate coordinates for $V_{\Phi}\cong\mathrm{T}^{n-1}$).

\vspace*{0.25in} \noindent {\bf Definition}: If $\{e_1,\ldots,e_n\}$
denotes the canonical basis of vectors in $\mathbb{R}^n$, we define
the vectors $v_1,\ldots,v_n$ by ${\displaystyle
v_1=\sum_{k=1}^n\,e_k}$, and for $j=2,\ldots,n$,
\[
v_j=v_1-\sum_{\ell=0}^{j-2}\,2\,e_{n-\ell}.
\]
By construction, the set $\{v_1,\ldots,v_n\}$ is linearly
independent, and so the $n\times n$ matrix ${\cal I}$, whose
$j^{\mbox{\small th}}$ column is the vector $v_j^T$, is invertible.
We also define the $n\times n$ invertible matrix ${\cal U}_j$ to be
the diagonal matrix whose $k^{\mbox{\small th}}$ diagonal element is
the $k^{\mbox{\small th}}$ component of the vector $v_j$ (in
particular, ${\cal U}_1$ is the identity matrix).  Note also that
${\cal U}_j^{-1}={\cal U}_j$, $j=1,\ldots,n$.

\vspace*{0.25in} Consider the following point in
$\mathrm{V}_{\Phi}\times\mathrm{V}_{\Psi}$:
\[
{\displaystyle
(\widehat{\Phi},\widehat{\Psi})=-\frac{\pi}{2}\,((v_1,\ldots,v_{n-1}),v_n)};
\]
then it is easy to compute that
\begin{equation}\label{matrix-I}
\widetilde{P}(\widehat{\Phi},\widehat{\Psi})=i\,{\cal I},
\end{equation}
where $\widetilde{P}$ is as in (\ref{real_suspended}).
If we define
\[
\widehat{A}^T\equiv (\hat{a}_1,\ldots,\hat{a}_n)^T={\cal
I}^{-1}\omega^T
\]
then
\[
F(\widehat{\Phi},\widehat{\Psi},\widehat{A};\omega)=0.
\]
Because the $\omega_j$ are rationally independent, it follows that
the components $\hat{a}_k$ of $\widehat{A}$ are all non-zero.

Furthermore, it is easy to compute the following derivatives
\[
J\equiv
D_{(\Psi,A)}F(\widehat{\Phi},\widehat{\Psi},\widehat{A};\omega)=
\left(\,\begin{array}{cc}\hat{a}_n\,{\cal U}_n&i\,{\cal I}\\[0.1in]\hat{a}_n\,{\cal U}_n&-i\,{\cal
      I}\end{array}\,\right),
\]
and
\[
K\equiv D_{\Phi}F(\widehat{\Phi},\widehat{\Psi},\widehat{A};\omega)=
\left(\,\begin{array}{cccc}\hat{a}_1\,{\cal U}_1&\hat{a}_2\,{\cal
      U}_2&\cdots&\hat{a}_{n-1}\,{\cal U}_{n-1}\\[0.1in]
\hat{a}_1\,{\cal U}_1&\hat{a}_2\,{\cal
  U}_2&\cdots&\hat{a}_{n-1}\,{\cal U}_{n-1}\end{array}\,\right).
\]

The matrix $J$ is invertible, since its inverse is easily computed
as
\[
J^{-1}=\left(\,\begin{array}{cc}\frac{1}{2\hat{a}_n}\,{\cal
      U}_{n}&\frac{1}{2\hat{a}_n}\,{\cal U}_n\\[0.1in]
-\frac{i}{2}\,{\cal I}^{-1}&\frac{i}{2}\,{\cal I}^{-1}
\end{array}\,\right).
\]
By the implicit function theorem, there exists a neighborhood ${N}$
of $\widehat{\Phi}$ in $\mathrm{V}_{\Phi}$ and a
unique smooth function
\[
\begin{array}{lll}
G:{N}&\longmapsto &\mathrm{V}_{\Psi}\times\mathbb{R}^n\\[0.15in]
G: \Phi&\longmapsto&G(\Phi)=(G_{\Psi}(\Phi),G_A(\Phi))
\end{array}
\]
such that
\[
G(\widehat{\Phi})=(\widehat{\Psi},\widehat{A})
\]
and
\begin{equation}
F(\Phi,G(\Phi);\omega)\equiv 0,\,\,\,\forall\,\Phi\in {N}.
\label{ift}
\end{equation}

Implicit differentiation of (\ref{ift}) yields that
\begin{equation}
\begin{array}{lll}
DG(\widehat{\Phi})&=&\left(\,\begin{array}{c}DG_{\Psi}(\widehat{\Phi})\\[0.1in]
DG_A(\widehat{\Phi})\end{array}\,\right)=
-J^{-1}K\\[0.25in]
&=&\left(\,\begin{array}{cccc} -\frac{\hat{a}_1}{\hat{a}_n}\,{\cal
U}_n\,{\cal U}_1& -\frac{\hat{a}_2}{\hat{a}_n}\,{\cal U}_n\,{\cal
U}_2&\cdots&
-\frac{\hat{a}_{n-1}}{\hat{a}_n}\,{\cal U}_n\,{\cal U}_{n-1}\\[0.1in]
\mbox{\bf 0}&\mbox{\bf 0}&\cdots&\mbox{\bf 0}\end{array}\,\right),
\end{array}
\label{diffeo}
\end{equation}
where $\mbox{\bf 0}$ denotes the $n\times n$ zero matrix.
Consequently,
\begin{equation}
DG_{\Psi}(\widehat{\Phi})=
\left(\,\begin{array}{cccc}-\frac{\hat{a}_1}{\hat{a}_n}\,{\cal
      U}_n\,{\cal U}_1
& -\frac{\hat{a}_2}{\hat{a}_n}\,{\cal U}_n\,{\cal U}_2&\cdots&
-\frac{\hat{a}_{n-1}}{\hat{a}_n}\,{\cal U}_n\,{\cal
  U}_{n-1}\end{array}\,\right),
\label{submerse}
\end{equation}
and it follows that the mapping
\[
G_{\Psi}:N\longrightarrow \mathrm{V}_{\Psi}
\]
is regular at $\widehat{\Phi}$.

Consider the following $n-1$ vectors in $(\,\mathbb{R}^n\,)^{n-1}$:
\[
\begin{array}{lll}
W_1&=&(\omega,0,\ldots,0,0)\\[0.1in]
W_2&=&(0,\omega,0,\ldots,0,0)\\[0.1in]
&\vdots&\\[0.1in]
W_{n-2}&=&(0,0,\ldots,\omega,0)\\[0.1in]
W_{n-1}&=&(0,0,\ldots,0,\omega),
\end{array}
\]
where $0$ represents the $0$ vector in $\mathbb{R}^n$, and we recall
that $\omega=(\omega_1,\ldots,\omega_n)$.

The set $\{W_1,\ldots,W_{n-1}\}$ is linearly independent, so for any
$\Phi\in\,\mathrm{V}_{\Phi}$, the set
\[
{\cal S}_{\Phi}=\{\,\Phi+\sum_{j=1}^{n-1}\sigma_j\,W_j\,\,(\mbox{\rm
  mod}\,\, \mathrm{V}_{\Phi})\,\,\,|\,\,\,0\leq |\sigma_j| < <
  1,\,\,j=1,\ldots,n\,\}
\]
is a small $n-1$-dimensional surface through $\Phi$ in
$\mathrm{V}_{\Phi}$.  We are interested in showing that for
$\Phi$ close enough to $\widehat{\Phi}$ in $N$, the image of ${\cal
  S}_{\Phi}$ by $G_{\Psi}$ in $\mathrm{V}_{\Psi}$ is transverse to the
integral curves of the vector field $\dot{\Psi}=\omega$. To show
this, we consider the function
\[
{\cal T}: N\longrightarrow \mathbb{R}
\]
defined by
\begin{equation}
{\cal T}(\Phi)=\mbox{\rm det}\left(\,DG_{\Psi}(\Phi)\cdot
  W_1^T\,\,\,\,\,\,\,DG_{\Psi}(\Phi)\cdot
  W_2^T\,\,\,\,\,\,\,\cdots\,\,\,\,\,\,\,DG_{\Psi}(\Phi)\cdot
  W_{n-1}^T\,\,\,\,\,\,\,\omega^T\,\right).
\label{Tdef}
\end{equation}
Obviously, ${\cal T}$ is continuous, and
\[
\begin{array}{lll}
{\cal T}(\widehat{\Phi})&=&\mbox{\rm det}
\left(\,\begin{array}{ccccc}-\frac{\hat{a}_1}{\hat{a}_n}\,{\cal
      U}_n\,{\cal U}_1\omega^T
& -\frac{\hat{a}_2}{\hat{a}_n}\,{\cal U}_n\,{\cal
U}_2\omega^T&\cdots& -\frac{\hat{a}_{n-1}}{\hat{a}_n}\,{\cal
U}_n\,{\cal
  U}_{n-1}\omega^T&{\cal U}_n\,{\cal U}_n\omega^T\end{array}\,\right)\\[0.15in]
&=&{\displaystyle\frac{(-1)^{n-1}}{\hat{a}_n^{n-1}}}\,\,\,\,\mbox{\rm
  det}\,{\cal U}_n\,\mbox{\rm det}
\left(\,\begin{array}{ccccc}\hat{a}_1\,{\cal U}_1\omega^T &
\hat{a}_2\,{\cal U}_2\omega^T&\cdots& \hat{a}_{n-1}\,{\cal
  U}_{n-1}\omega^T&{\cal
  U}_n\omega^T\end{array}\,\right)\\[0.15in]
&=&{\displaystyle\frac{(\omega_1\,\omega_2\,\cdots\,\omega_n)\,(\hat{a}_1\,\hat{a}_2\,\cdots\,\hat{a}_{n-1})}{\hat{a}_{n}^{n-1}}}\,\,\mbox{\rm
    det}\,{\cal I}\\[0.15in]
&\neq&0.
\end{array}
\]
It follows that there is a neighborhood $N'\subseteq N$ in which
${\cal T}\neq 0$.  This is equivalent to saying that for all
$\Phi\in N'$, the image of ${\cal
  S}_{\Phi}$ by $G_{\Psi}$ in $\mathrm{V}_{\Psi}$ is transverse to the
integral curves of the vector field $\dot{\Psi}=\omega$.

For each $j=1,\ldots,n-1$, the integral curves of the vector field
$\dot{\Phi}^j=\omega$ are dense. Thus,
for any $\varepsilon>0$, there is a $\tau_{j,\varepsilon}>0$ and an
$s_{j,\varepsilon}>0$ such that the integral curve segment
\[
\{\,\Phi^j=\tau_j\omega\,\,\,(\mbox{\rm
mod}\,\,(\mathbb{S}^1)^n)\,\,\,|\,\,\,\tau_{j,\varepsilon}-s_{j,\varepsilon}<\tau_j<\tau_{j,\varepsilon}+s_{j,\varepsilon}\,\}
\]
is in the $\varepsilon$-ball centered on $\frac{\pi}{2}v_j$ in
$(\mathbb{S}^1)^n$. For $\varepsilon>0$ small enough, the surface
\[
\{\,\Phi=(\tau_1\omega,\tau_2\omega,\ldots,\tau_{n-1}\omega)\,\,(\mbox{\rm
  mod}\,\,\mathrm{V}_{\Phi})\,\,\,|\,\,\,\tau_{j,\varepsilon}-s_{j,\varepsilon}<\tau_j<\tau_{j,\varepsilon}+s_{j,\varepsilon}\,\}
\]
is contained in $N'$ and coincides with the surface ${\cal
  S}_{\Phi^*}$ for
\[
\Phi^*=(\tau_{1,\varepsilon}\omega,\tau_{2,\varepsilon}\omega,\ldots,\tau_{(n-1),\varepsilon}\omega)\,\,\,\,(\mbox{\rm
  mod}\,\,\mathrm{V}_{\Phi}).
\]
  Thus, by our previous
  result,
the $n-1$-dimensional
  surface $G_{\Psi}({\cal S}_{\Phi^*})$ is transverse to the integral
  curves of $\dot{\Psi}=\omega$ in $\mathrm{V}_{\Psi}$.  Since these
  integral curves are dense in $\mathrm{V}_{\Psi}$, there are infinitely many intersections
  with $G_{\Psi}({\cal S}_{\Phi^*})$ near
  the point $\widehat{\Psi}=G_{\Psi}(\widehat{\Phi})$.

Let $\stackrel{\circ}{\Psi}\in G_{\Psi}({\cal S}_{\Phi^{*}})$ be
such an intersection point near $\widehat{\Psi}$.  Then there is a
$\stackrel{\circ}{\tau}_n>0$ such that
\[
\stackrel{\circ}{\Psi}=\stackrel{\circ}{\tau_n}\omega\,\,(\mbox{\rm
  mod}\,\,\mathrm{V}_{\Psi}).
\]
Let $\stackrel{\circ}{\Phi}\in {\cal S}_{\Phi^*}$ be such that
$G_{\Psi}(\stackrel{\circ}{\Phi})=\stackrel{\circ}{\Psi}$. Then
there are
$\stackrel{\circ}{\tau}_1>0,\,\stackrel{\circ}{\tau}_2>0,\,\ldots\,,\stackrel{\circ}{\tau}_{n-1}>0$
such that
\[
\stackrel{\circ}{\Phi}=(\stackrel{\circ}{\tau}_1\omega,\stackrel{\circ}{\tau}_2\omega,\ldots,\stackrel{\circ}{\tau}_{n-1}\omega)\,\,(\mbox{\rm
  mod}\,\,\,\mathrm{V}_{\Phi}).
\]

It follows from (\ref{ift}) that
$F(\stackrel{\circ}{\Phi},\stackrel{\circ}{\Psi},G_A(\stackrel{\circ}{\Phi});\omega)=0$,
and by construction, this corresponds to a solution of
(\ref{chareq_matrix}). \qed

The next theorem shows that the previous realization result holds
for open sets near solutions found in Theorem~\ref{thm:main}.

\begin{theorem}\label{thm:open} Suppose
$\omega_1>0,\,\omega_2>0,\,\ldots\,,\omega_{n}>0$ are linearly
independent over the rationals.  There exists a neigborhood ${\cal N}$ of
$\omega=(\omega_1,\ldots,\omega_n)$ in $\mathbb{R}^n$ and a smooth
mapping
\[\begin{array}{llll}
H:&V&\longrightarrow&\mathbb{R}^n\times\mathbb{R}^n\\[0.15in]
&\omega&\longmapsto&H(\omega)=(\tau(\omega),A(\omega))=((\tau_1(\omega),\ldots,\tau_n(\omega)),(a_1(\omega),\ldots,a_n(\omega)))
\end{array}
\]
such that
\begin{equation}
\begin{array}{ccc}
{\displaystyle\sum_{k=1}^n\,a_k(\omega)\,e^{-i\omega_j\tau_k(\omega)}}&=&i\omega_j,\,\,\,\,\,\,j=1,\ldots,n\\[0.15in]
{\displaystyle\sum_{k=1}^n\,a_k(\omega)\,e^{i\omega_j\tau_k(\omega)}}&=&-i\omega_j,\,\,\,\,\,\,j=1,\ldots,n
\end{array}
\label{chareq_ift}
\end{equation}
for all $\omega\in {\cal N}$.
\end{theorem}

\proof We consider the system $F=0$ given by (\ref{real_suspended}).
We have already shown in Theorem~\ref{thm:main} that, for fixed
$\omega$ linearly independent over the rationals, there exists
infinitely many solutions to $F=0$.  We again use an implicit
function theorem argument combined with the density of irrational
torus flows.

Consider the mapping
\begin{equation}
\begin{array}{llll}
Q&:\mathbb{R}^n\times\mathbb{R}^n\times\mathbb{R}^n&\longrightarrow&\mathbb{R}^n\times\mathbb{R}^n\\[0.15in]
&(\tau,A,\omega)&\longmapsto&Q(\tau,A,\omega)=F((\tau_1\omega,\ldots,\tau_{n-1}\omega),\tau_n\omega,A;\omega),
\end{array}
\label{main_sys}
\end{equation}
where $F$ is as in (\ref{real_suspended}). Therefore
\begin{equation}
D_{\tau}Q(\tau,A,\omega)=D_{((\Phi^1,\ldots,\Phi^{n-1}),\Psi)}F((\tau_1\omega,\ldots,\tau_{n-1}\omega),\tau_n\omega,A;\omega)\cdot
\left(\begin{array}{ccccc}\omega^T&0&0&\cdots&0\\[0.15in]
0&\omega^T&0&\cdots&0\\[0.15in]
\vdots&\vdots&\vdots&\vdots&\vdots\\[0.15in]
0&0&0&\cdots&\omega^T
\end{array}
\right), \label{der1}
\end{equation}
where each $0$ in the matrix above is an $n$-dimensional zero column
vector; and
\begin{equation}
D_{A}Q(\tau,A,\omega)=D_AF((\tau_1\omega,\ldots,\tau_{n-1}\omega),\tau_n\omega,A;\omega).
\label{der2}
\end{equation}
Thus, we wish to show that the $2n\times 2n$ matrix
\begin{equation}
\left(\,D_{\tau}Q(\tau,A,\omega)\,\,\,\,\,\,\,\,\,\,\,D_{A}Q(\tau,A,\omega)\,\right)
\label{der3}
\end{equation}
is invertible at the solutions to (\ref{chareq}) we have found in
Theorem~\ref{thm:main}.

For positive integers $p$ and $q$, let $\mbox{\rm Mat}_{p,q}$ denote
the space of $p\times q$ matrices.  Consider the following mappings
associated to (\ref{der1}), (\ref{der2}) and (\ref{der3}):
\[
{\cal
  R}_1:\mathrm{V}_{\Phi}\times\mathrm{V}_{\Psi}\times\mathbb{R}^n\times\mathbb{R}^n\longrightarrow\mbox{\rm
  Mat}_{2n,n}
\]
defined by
\[
{\cal R}_1(\Phi,\Psi,A,\omega)=
D_{(\Phi,\Psi)}F(\Phi,\Psi,A;\omega)\cdot
\left(\begin{array}{ccccc}\omega^T&0&0&\cdots&0\\[0.15in]
0&\omega^T&0&\cdots&0\\[0.15in]
\vdots&\vdots&\vdots&\vdots&\vdots\\[0.15in]
0&0&0&\cdots&\omega^T
\end{array}
\right),
\]
\[
{\cal
  R}_2:\mathrm{V}_{\Psi}\times\mathrm{V}_{\Psi}\times\mathbb{R}^n\times\mathbb{R}^n\longrightarrow\mbox{\rm
  Mat}_{2n,n}
\]
defined by
\[
{\cal
  R}_2(\Phi,\Psi,A,\omega)=D_{A}F(\Phi,\Psi,A;\omega),
\]
and
\[
{\cal
  R}:\mathrm{V}_{\Phi}\times\mathrm{V}_{\Psi}\times\mathbb{R}^n\times\mathbb{R}^n\longrightarrow\mbox{\rm
  Mat}_{2n,2n}
\]
defined by
\[
{\cal R}(\Phi,\Psi,A,\omega)=\left( \,{\cal
  R}_1(\Phi,\Psi,A,\omega)\,\,\,\,\,\,\,\,\,\,\,\,\,\,\,
{\cal R}_2(\Phi,\Psi,A,\omega) \,\right).
\]
Now, a simple computation (similar to those done in the proof of
Theorem~\ref{thm:main}) shows that
\[
{\cal
  R}\left(-\frac{\pi}{2}(v_1,\ldots,v_n),A,\omega\right)=\left(\begin{array}{cc}
{\cal Z}&i\,{\cal I}\\[0.15in]
{\cal Z}&-i\,{\cal I}
\end{array}\right),
\]
where
\[
{\cal Z}=\left(\begin{array}{llll} a_1\,{\cal
U}_1\omega^T&a_2\,{\cal U}_2\omega^T&\cdots&a_n\,{\cal
  U}_n\omega^T\end{array}\right).
\]
If none of the $a_j$ vanish, then the $n\times n$ matrix ${\cal Z}$
is invertible, since its determinant is
\[
\mbox{\rm det}\,{\cal Z}=\prod_{j=1}^n\,a_j\omega_j\,\,\mbox{\rm
  det}\,{\cal I}\neq 0.
\]
Thus,
\[
{\cal
  R}\left(-\frac{\pi}{2}(v_1,\ldots,v_n),A,\omega\right)^{-1}=\left(\begin{array}{cc}
\frac{1}{2}{\cal Z}^{-1}&\frac{1}{2}{\cal Z}^{-1}\\[0.15in]
-\frac{i}{2}\,{\cal I}^{-1}&\frac{i}{2}\,{\cal I}^{-1}
\end{array}\right).
\]
By continuity, there is thus a neighborhood ${\cal N}$ of the point
${\displaystyle
  -\frac{\pi}{2}(v_1,\ldots,v_n)}$ in $\mathrm{V}_{\Phi}\times\mathrm{V}_{\Psi}$ in which
${\cal R}$ is invertible.  By Theorem~\ref{thm:main}, there are
infinitely many solutions of $Q=0$ (see (\ref{main_sys})) in ${\cal
N}$, and the Jacobian matrix (\ref{der3}) is thus invertible at
these solutions. We get the conclusion of Theorem~\ref{thm:open} by
the implicit function theorem.\qed

\subsection{Example: $\D_3$-symmetric system}

Theorem~\ref{thm:main} is written in the context of scalar
delay-differential equations. However, in this section, we look at
an example of a $\D_3$-symmetric system of delay-differential
equations where Theorem~\ref{thm:main} can be applied and then
proceed to explain the generalization of this theorem which has
applications to symmetric systems of delay-differential equations.

\begin{example}
Let $\Gamma=\D_3$ be the group generated by $\kappa$ and $\gamma$
act on $\R^{3}$ as follows:
\[
\kappa.(x_1,x_2,x_3)=(x_1,x_3,x_2),\qquad
\gamma.(x_1,x_2,x_3)=(x_3,x_1,x_2).
\]
Consider a linear $\D_3$-symmetric coupled cell system with delayed
coupling where each cell is one-dimensional and has the following
form.
\begin{equation}\label{eq:D3}
\begin{array}{rcl}
\dot x_1&=&a_1 x_1(t-\tau_1)+a_2[x_2(t-\tau_2)+x_3(t-\tau_2)] \\
\dot x_2&=&a_1 x_2(t-\tau_1)+a_2[x_3(t-\tau_2)+x_1(t-\tau_2)] \\
\dot x_3&=&a_1 x_3(t-\tau_1)+a_2[x_1(t-\tau_2)+x_2(t-\tau_2)].
\end{array}
\end{equation}
where $x_i\in \R$ for $i=1,2,3$ and $a_1,a_2,a_3\in \R$. The
characteristic equation of system~(\ref{eq:D3}) is obtained by
substituting $(x_1,x_2,x_3)=(w_1 e^{\lambda t},w_2 e^{\lambda t},w_3
e^{\lambda t})$ into the equations. We obtain after simplification
\[
\begin{array}{rcl}
\lambda w_1&=& a_1 e^{-\lambda \tau_1}w_1+a_2e^{-\lambda \tau_2}[w_2+w_3] \\
\lambda w_2&=& a_1 e^{-\lambda \tau_1}w_2+a_2e^{-\lambda \tau_2}[w_3+w_1]\\
\lambda w_3&=& a_1 e^{-\lambda \tau_1}w_3+a_2e^{-\lambda
\tau_2}[w_1+w_2]
\end{array}
\]
and rearranging the terms we have
\begin{equation}\label{eq:D3char1}
\left[(\lambda-a_1 e^{-\lambda \tau_1})I-a_2
e^{-\lambda\tau_2}\left(\begin{array}{ccc} 0 & 1 & 1 \\ 1 & 0 & 1 \\
1 & 1 & 0
\end{array}\right)\right]\left(\begin{array}{c} w_1\\ w_2 \\
w_3\end{array}\right)=0.
\end{equation}
where $I$ is the $3\times 3$ identity matrix. Letting
$\alpha=\lambda-a_1 e^{-\lambda \tau_1}$ and $\beta=-a_2
e^{-\lambda\tau_2}$ equation~(\ref{eq:D3char1}) becomes
\[
\left(\begin{array}{ccc} \alpha & \beta & \beta \\ \beta & \alpha & \beta \\
\beta & \beta & \alpha
\end{array}\right)
\left(\begin{array}{c} w_1\\ w_2 \\
w_3\end{array}\right)=0.
\]
Let
\[
\Delta(\lambda)=\left(\begin{array}{ccc} \alpha & \beta & \beta \\ \beta & \alpha & \beta \\
\beta & \beta & \alpha
\end{array}\right).
\]
We complexify $\R^{3}$ and look at the isotypic decomposition of
$\C^{3}$ by the action of $\D_3$:
\[
\C^{3}=V_0\oplus V_1\oplus V_2
\]
where $V_0$ is the trivial representation of $\D_3$ and $V_1$, $V_2$
are the standard irreducible representations of $\D_3$ (all
representations are one-dimensional complex). A basis for $V_0$ is
$u_0=(v,v,v)^{t}$, a basis for $V_1$ is $u_1=(v,e^{2\pi
i/3}v,e^{4\pi i/3}v)$ and a basis for $V_2$ is $u_2=(v,e^{4\pi
i/3}v,e^{2\pi i/3}v)$. Therefore,
\[
\Delta(\lambda)u_0=(\alpha+2\beta) u_0
\]
and
\[
\Delta(\lambda)u_1=(\alpha-\beta) u_1,\qquad
\Delta(\lambda)u_2=(\alpha-\beta) u_2
\]
since $e^{4\pi i/3}=\ov{e^{2\pi i/3}}$. Therefore, in the basis
given by the isotypic decomposition of $\C^{3}$, $\Delta(\lambda)$
block diagonalizes so that we have
\[
\left(\begin{array}{ccc}\alpha+2\beta & 0 & 0 \\ 0 & \alpha-\beta & 0 \\
0 & 0 & \alpha-\beta\end{array}\right)\left(\begin{array}{c} \tilde{w}_1\\ \tilde{w}_2 \\
\tilde{w}_3\end{array}\right)=0.
\]
Hence, the eigenvalues are solutions to
\[
\det \Delta(\lambda)=(\alpha+2\beta)(\alpha-\beta)^2=(\lambda-a_1
e^{-\lambda \tau_1}-2a_2 e^{-\lambda\tau_2})(\lambda-a_1 e^{-\lambda
\tau_1}+a_2 e^{-\lambda\tau_2})^2=0.
\]
Each factor of the characteristic equation is of the same form as
the characteristic equation for a scalar delay-differential
equation. Therefore, by letting $\tilde{a}_1=a_1$ and
$\tilde{a}_2=2a_2$ in $(\lambda-a_1 e^{-\lambda \tau_1}-2a_2
e^{-\lambda\tau_2})$, Theorem~\ref{thm:main} applies directly. The
same is true for the factor $(\lambda-a_1 e^{-\lambda \tau_1}+a_2
e^{-\lambda\tau_2})$ where we let $\tilde{a}_1=a_1$ and
$\tilde{a}_2=-a_2$. Hence, for any choice of a set of complex
numbers $\Lambda=\{i\omega_1,i\omega_2\}$ with $\omega_1,\omega_2>0$
and rationally independent, there exists a linear $\D_3$ symmetric
coupled cell system including $\Lambda$ in its spectrum.
\end{example}

In the context of bifurcation theory, the symmetry properties of the
critical eigenspace depends on which factor contains the critical
eigenvalue and this leads to different bifurcation behaviour. Two
imaginary eigenvalues in the first factor corresponds to a
nonresonant Hopf/Hopf mode interaction (without symmetry) while the
second case leads to a nonresonant $\D_3$ Hopf/Hopf mode
interaction. Details of the unfolding of these bifurcations can be
found respectively in Kuznetsov~\cite{Kuznetsov} and Golubitsky~{\em
et. al.}~\cite{GSS88}.

Note that Theorem~\ref{thm:main} is not sufficient to guarantee the
existence of a linear $\D_3$ symmetric coupled cell system with
$i\omega_1$ satisfying the first factor and $i\omega_2$ satisfying
the second factor simultaneously. We characterize this situation as
follows. Let $b_1^1=b_1^2=1$ and $b_2^1=2$ and $b_2^2=-1$ and for
fixed rationally independent $i\omega_1,i\omega_2$ (with
$\omega_1,\omega_2>0$), we look for $a_1,a_2$ and $\tau_1,\tau_2$
such that
\begin{equation}\label{ex:char-eq}
\begin{array}{l}
a_1 b_1^1 e^{-i\omega_1\tau_1}+a_2 b_2^1 e^{-i\omega_1\tau_2}=i\omega_1\\
a_1 b_1^2 e^{-i\omega_2 \tau_1}+a_2 b_2^2
e^{-i\omega_2\tau_2}=i\omega_2
\end{array}
\end{equation}
and their complex conjugate equations are satisfied. This is the
context of the next theorem which is a generalization of
Theorem~\ref{thm:main}. We state this result in a general form below
and postpone the proof to Section~\ref{sec:main-proof} as it follows
similar steps as the proof of Theorem~\ref{thm:main}.

Note that in the proof of Theorem~\ref{thm:main}, the matrix ${\cal
I}$ defined in~(\ref{matrix-I}) is nonsingular by construction and
this is a crucial step in the argument. For this more general result
we shall present, the matrix which holds a similar role is denoted
by ${\cal I}_{B}$ since it is a matrix consisting of $\pm$ the
constants $b_{k}^{j}$ which appear in equations~(\ref{ex:char-eq}).
The form of this matrix is not relevant for the moment and the
structure of the matrix is described in
Section~\ref{sec:main-proof}. We are now ready to state the theorem.

\begin{theorem}\label{thm:main2}
Consider the factors
\begin{equation}\label{char-eq:product}
\prod_{j=1}^{r} \left(\lambda-\sum_{k=1}^{n} a_{k} b_{k}^{j}
e^{\lambda\tau_{k}}\right)
\end{equation}
of a characteristic polynomial where the constants $b_k^{j}\in
\R\setminus\{0\}$ are fixed for all $j=1,\ldots,r$, $k=1,\ldots,n$
and suppose that $\det {\cal I}_{B}\neq 0$. Suppose that
$\omega_1>0,\,\omega_2>0,\,\ldots\,,\omega_{n}>0$ are linearly
independent over the rationals. Then there exists
$\tau_1>0,\,\tau_2>0,\,\ldots\,,\tau_n>0$,
$a_1\in\mathbb{R},\,a_2\in\mathbb{R},\,\ldots\,,a_n\in\mathbb{R}$
such that for all $j=1,\ldots,r$,
\[
\left(\lambda-\sum_{k=1}^{n} a_{k}b_k^{j}
e^{\lambda\tau_{k}}\right)=0
\]
has roots $i\omega_{\ell}^{j}$ for $\ell=1,\ldots,\ell_{j}$ where
$\ell_1+\cdots+\ell_{r}=n$.
\end{theorem}

This theorem is applied in the following sections to the case of
$\D_n$ symmetric coupled one-dimensional cell systems. If $n$ odd,
it is easy to show that $b_{k}^{j}\neq 0$ holds, but for $n$ even,
several of the $b_{k}^{j}$'s can be zero and Theorem~\ref{thm:main2}
cannot be applied directly.

\section{Linear $\Gamma$-symmetric delay-differential equations}
\label{sec:symmetry} For the results of this section, we find it
convenient to introduce the well know abstract setting, see for
instance Hale and Lunel~\cite{HL}, adapted to the symmetric case.
Let $C_n=C([-\tau,0],{\mathbb C}^n)$ be the Banach space of
continuous functions from the interval $[-\tau,0]$, into ${\mathbb
C}^n$ ($\tau>0$) endowed with the norm of uniform convergence.
Consider the linear homogeneous RFDE
\begin{equation}
\dot{z}(t)={\cal L}_0(z_t), \label{fde}
\end{equation}
where ${\cal L}_0$ is a bounded linear operator from $C_n$ into
${\mathbb C}^n$. We write
\[
{\cal L}_0(\varphi)=\int_{-\tau}^{0}\,d\eta(\theta)\varphi(\theta),
\]
where $\eta$ is an $n\times n$ matrix-valued function of bounded
variation defined on $[-\tau,0]$. The characteristic equation is
\begin{equation}
\mbox{\rm det}\,\Delta(\lambda)=0,\;\;\;\;\;\;\;\; \mbox{\rm
where}\;\;\Delta(\lambda)=\lambda\,I_n-\int_{-\tau}^0\,
d\eta(\theta)e^{\lambda\theta}, \label{char_eq}
\end{equation}
where $I_n$ is the $n\times n$ identity matrix. Note that
$e^{\lambda\theta}=e^{\lambda\theta}I_n$.

Suppose that $\Gamma$ is a compact group of transformations acting
linearly on ${\mathbb C}^n$.  We say that (\ref{fde}) is {\em
$\Gamma$-equivariant} if
\begin{equation}
\gamma\cdot \eta(\theta)=\eta(\theta)\cdot\gamma,\,\,\,\forall\,
\gamma\in\Gamma, \,\theta\in [-\tau,0]. \label{equiv_cond}
\end{equation}
The group action of $\Gamma$ on $\C^{n}$ induces an isotypic
decomposition of $\C^{n}$:
\[
\C^{n}=V_1\oplus V_2\oplus \cdots V_k
\]
where $V_i=U_i\oplus \cdots U_i$ for irreducible representations
$U_i$ of $\Gamma$ and $U_i\not\simeq U_j$ for $i\neq j$. Since
$\eta(\theta)$ commutes with the action of $\Gamma$, then
\[
\eta(\theta)V_i\subset V_i
\]
for all $i=1,\ldots,k$.

Therefore, $\Delta(\lambda)$ also commutes with the representation
of $\Gamma$. Indeed, for all $\gamma\in \Gamma$
\[
\begin{array}{rcl}
\Delta(\lambda)\gamma&=&\lambda I\gamma-\left[\int_{-\tau}^{0}
d\eta(\theta)e^{\lambda\theta}\right]\gamma\\
&=&\gamma \lambda I-\left[\int_{-\tau}^{0} d\eta(\theta)\gamma
e^{\lambda\theta}\right]\\
&=&\gamma \lambda I-\left[\int_{-\tau}^{0} \gamma d\eta(\theta)
e^{\lambda\theta}\right]\\
&=&\gamma \left(\lambda I-\int_{-\tau}^{0} d\eta(\theta)
e^{\lambda\theta}\right)=\gamma \Delta(\lambda).
\end{array}
\]
Thus,
\[
\Delta(\lambda)V_i\subset V_i
\]
and for all $i=1,\ldots,k$ and we can write $\Delta(\lambda)$ in
block diagonal form:
\[
\Delta(\lambda)=\diag(\Delta_1(\lambda),\ldots,\Delta_k(\lambda)).
\]
The characteristic equation then becomes
\[
\det\Delta(\lambda)=\prod_{i=1}^{k} \det\Delta_{i}(\lambda).
\]
Therefore we are led to the following result.
\begin{proposition}\label{prop:1d}
Suppose that $V_i=U_i$ and $U_i$ is a one-dimensional irreducible
representation of $\Gamma$. Then
\[
\det\Delta_{i}(\lambda)=\lambda-\sum_{j=1}^{\ell}
a_{j}e^{\lambda\tau_j}
\]
\end{proposition}

\begin{corollary}
Theorem~\ref{thm:main} applies to factors of the characteristic
equation which correspond to the context of
Proposition~\ref{prop:1d}.
\end{corollary}

\subsection{Delayed coupled cell systems with
$\D_n$-symmetry}\label{section:delay-ring} Multiple authors have
studied Hopf bifurcation in $\D_n$ symmetric rings of cells with
delayed coupling where each cell is one-dimensional. The
differential equation systems in those papers have the following
general form. For $i=1,\ldots,n$, the dynamics of cell $i$ is given
respectively for $n$ odd and $n$ even by:
\begin{equation}\label{eq:Dn-syst-odd}
\dot
x_{i}(t)=f(X_i)+g(x_{i+1},\ldots,x_{i+(n-1)/2},x_{i-(n-1)/2},\ldots,x_{i-1})
\end{equation}
\begin{equation}\label{eq:Dn-syst-even}
\dot
x_{i}(t)=f(X_i)+g(x_{i+1},\ldots,x_{i+(n/2-1)},x_{i+n/2},x_{i-(n/2-1)},\ldots,x_{i-1})
\end{equation}
where $X_i=(x_{i}(t-s_{1}),\ldots,x_{i}(t-s_{m}))$,
$x_{j}=x_{j}(t-\tau_j)$ for $j\neq i$, $f:\R^{m}\to\R$,
$g:\R^{n-1}\to\R$ are smooth functions and $\tau_j,s_{\ell}\in
[0,\tau]$ for all $j\neq i$ and $s=1,\ldots,m$. Here, $f$ is called
the internal dynamics and $g$ is the coupling function. In the
papers described above, these systems are $\D_n$-equivariant.

\subsubsection{Characterization of delayed $\D_n$ networks}
We introduce a more general notation for delayed symmetrically
coupled cell systems inspired by recent work on (non-necessarily
symmetric) coupled cell systems of ordinary differential equations,
see for instance~\cite{Leite-Golubitsky06}. Suppose that each cell
in the system has phase space $\R^{k}$. We generalize
systems~(\ref{eq:Dn-syst-even}) and~(\ref{eq:Dn-syst-odd}) to
\begin{equation}\label{eq:Dn-syst-odd1}
\dot
X_{i}(t)=f(\widetilde{X}_i,\widetilde{X}_{i+1},\ldots,\widetilde{X}_{i+(n-1)/2},\widetilde{X}_{i-(n-1)/2},\ldots,\widetilde{X}_{i-1}),
\qquad i=1,\ldots,n
\end{equation}
and
\begin{equation}\label{eq:Dn-syst-even1}
\dot
X_{i}(t)=f(\widetilde{X}_i,\widetilde{X}_{i+1},\ldots,\widetilde{X}_{i+(n/2-1)},\widetilde{X}_{i+n/2},\widetilde{X}_{i-(n/2-1)},\ldots,
\widetilde{X}_{i-1}), \qquad i=1,\ldots,n
\end{equation}
where
\[
\widetilde{X}_{j}=(X_j(t-\tau_1),\ldots,X_j(t-\tau_m)),
\]
$f:(\R^{m})^{n}\to \R^{k}$ is a smooth function and the position of
$\widetilde{X}_{k}$ corresponds to the coupling from cell $k$ to
cell $i$. We say that cells $j$ and $k$ have {\em identical
coupling} to cell $i$ if
\[
f(X_i,\ldots,u,\ldots,v,\ldots)=f(X_i,\ldots,v,\dots,u,\ldots)
\]
where $u$ and $v$ are permuted from positions $j$ and $k$. We
rewrite systems~(\ref{eq:Dn-syst-odd1}) and~(\ref{eq:Dn-syst-even1})
as
\begin{equation}\label{eq:total}
\dot X=F(\widetilde{X})
\end{equation}
where $X=(X_1,\ldots,X_n)^{t}$
\[
\begin{array}{rcl}
\widetilde{X}&=&\widetilde{X}_i,\widetilde{X}_{i+1},\ldots,\widetilde{X}_{i+(n-1)/2},
\widetilde{X}_{i-(n-1)/2},\ldots,\widetilde{X}_{i-1},\quad \mbox{or}\\
\vspace{1\baselineskip}
\widetilde{X}&=&\widetilde{X}_i,\widetilde{X}_{i+1},\ldots,\widetilde{X}_{i+(n/2-1)},
\widetilde{X}_{i+n/2},\widetilde{X}_{i-(n/2-1)},\ldots,
\widetilde{X}_{i-1},
\end{array}
\]
and $F:(\R^{m})^{n}\to \R^{nk}$ has $i^{th}$ component given by the
formulas above for $\dot X_i(t)$ (for either $n$ odd or even).

\begin{figure}[ht]
\centerline{%
\epsfxsize=0.4\hsize\epsffile{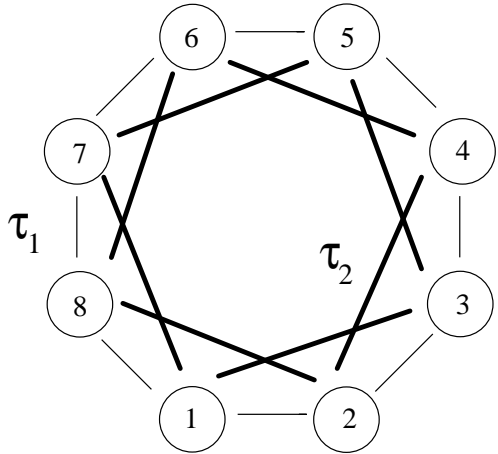}} \label{graph:Dn}
\caption{Representation of a $\D_8$-symmetric coupled cell system
with nearest neighbour with delay $\tau_1$ and second nearest
neighbour coupling with delay $\tau_2$.}
\end{figure}

Consider the group $\D_n$, with generators $\rho$ and $\kappa$
acting on $\R^{kn}$ as follows:
\begin{equation}\label{Dn-action}
\begin{array}{ll}
\rho.(X_1,\ldots,X_n)=(X_n,X_1,X_2,\ldots,X_{n-1}) & \\
\kappa.(X_1,\ldots,X_n)=(X_1,X_n,\ldots,X_{n+2-j},\ldots,X_{(n+1)/2},X_{(n+3)/2},\ldots,X_{j},\ldots,X_{2})
& \mbox{if $n$ is odd} \\
\kappa.(X_1,\ldots,X_n)=(X_n,\ldots,X_{n+1-j},\ldots,X_{n/2},X_{n/2+1},\ldots,X_j,\ldots,X_{1})&
\mbox{if $n$ is even}.
\end{array}
\end{equation}
These networks can be represented by a graph as shown in
Figure~\ref{graph:Dn}. In this case, we have an eight-cell network
with bidirectional nearest and next-nearest neighbour coupling with
delays $\tau_1$ and $\tau_2$ respectively.

Without loss of generality we assume that the networks are
transitive. That is, all cells in the network can be reached from
any other cell via the coupling arrows.

We now characterize the connections in the network so that the
delay-differential system is $\D_n$-symmetric. We think of each cell
in the network has having $[(n-1)/2]$ neighbours on each side and an
opposite cell if $n$ is even. Graphically, it is clear that a
$n$-cell network is $\D_n$ symmetric if for all cells in the
network, all connections to and from the $j^{th}$ neighbour on each
side (or the opposite cell if $n$ is even) are all the same; that
is, the coupling term and its delay must be the same for all those
connections. This idea is formalized in the next result.

\begin{proposition}\label{prop:Dn-net}
A transitive network of $n$ coupled identical cells, with delays, is
$\D_n$-equivariant if and only if it satisfies the conditions below.
\begin{enumerate}[(i)]
\item Suppose that cell $1$ receives an input from cell $j$ with delay
$\epsilon\in [0,\tau]$, then every cell $i$ in the network
($i=2,\ldots,n$) receives an input from cell $(i+j-1)\,\mod n$ with
delay $\epsilon$ identical to the one received by cell $1$.

\item For every connection in part (i), there is an identical
connection from cell $i$ to cell $(i+j-1)\,\mod n$ with delay
$\epsilon$.
\end{enumerate}
\end{proposition}
From Proposition~\ref{prop:Dn-net}, a transitive $\D_n$ symmetric
network with an even number of cells must have nearest neighbour
coupling.

\proof We begin by looking at $\rho$-equivariance. Denote by
$[w]_{i}$ the $i^{th}$ row of vector $w$. Then,
\[
[\rho
F(\widetilde{X})]_{i}=f(\widetilde{X}_{i-1},\widetilde{X}_i,\ldots,\widetilde{X}_{i-2}).
\]
and since $\rho X=(X_n,X_1,\ldots,X_{n-1})$ we have
\[
[F(\rho
\widetilde{X})]_{i}=f(\widetilde{X}_{i-1},\widetilde{X}_i,\ldots,\widetilde{X}_{i-2}).
\]
Thus, $\rho$-equivariance holds automatically by the structure of
the equations.

Consider the case $n$ odd (the case $n$ even has a similar proof and
is omitted). Suppose that cell $1$ receives an input from cell $j$.
We look at the system of equation~(\ref{eq:Dn-syst-odd1}) and focus
on the possible coupling from cell $(i+j-1)\,\mod n$ to cell $i$.
Moreover, consider the possible connection from cell $i$ to cell
$(i-j+1)\,\mod n$. Note that the connections from $(i+j-1)\,\mod n$
to $i$ and from $i$ to $(i-j+1)\,\mod n$ are obtained by taking the
index and subtracting $j-1$. Finally, consider the possible
connection from cell $(i-j+1)\,\mod n$ to cell $i$. We now show that
$F(\widetilde{X})$ is $\kappa$-equivariant (and so
$\D_n$-equivariant) if and only if the connections defined above are
identical. We see that
\[
\begin{array}{rcl}
\kappa F(\widetilde{X})&=&\kappa \left(\begin{array}{c}
f(\widetilde{X}_1,\ldots,\widetilde{X}_{j},\ldots,\widetilde{X}_{n+2-j},\ldots)\\
\vdots\\
f(\widetilde{X}_i,\ldots,\widetilde{X}_{i+j-1},\ldots,\widetilde{X}_{i-j+1},\ldots)\\
\vdots\\
f(\widetilde{X}_{n+2-i},\ldots,\widetilde{X}_{n+2-i+j-1},\ldots,\widetilde{X}_{n+2-i-j+1},\ldots)\\
\vdots
\end{array}\right)\\
\\
&=&\left(\begin{array}{c}
f(\widetilde{X}_1,\ldots,\widetilde{X}_{j},\ldots,\widetilde{X}_{n+2-j},\ldots)\\
\vdots\\
f(\widetilde{X}_{n+2-i},\ldots,\widetilde{X}_{n+2-i+j-1},\ldots,\widetilde{X}_{n+2-i-j+1},\ldots) \\
\vdots\\
f(\widetilde{X}_i,\ldots,\widetilde{X}_{i+j-1},\ldots,\widetilde{X}_{i-j+1},\ldots)\\
\vdots
\end{array}\right)
\end{array}
\]
and
\[
F(\kappa \widetilde{X})=\left(\begin{array}{c}
f(\widetilde{X}_1,\ldots,\widetilde{X}_{n+2-j},\ldots,\widetilde{X}_{j},\ldots)\\
\vdots\\
f(\widetilde{X}_{n+2-i},\ldots,\widetilde{X}_{n+2-(i+j-1)},\ldots,\widetilde{X}_{n+2-(i-j+1)},\ldots)\\
\vdots\\
f(\widetilde{X}_{i},\ldots,\widetilde{X}_{i-j+1},\ldots,\widetilde{X}_{i+j-1},\ldots)\\
\vdots
\end{array}\right).
\]
We now show that parts (i) and (ii) implies $\kappa$-equivariance.
If part (i) holds, the coupling from cell $(i+j-1)\,\mod n$ to cell
$i$ and the coupling from cell $i$ to cell $(i-j+1)\,\mod n$ are
identical. Then, by part (ii), the coupling from cell $(i-j+1)\,\mod
n$ to cell $i$ is identical to the coupling from cell $i$ to cell
$(i-j+1)\,\mod n$. Therefore, the coupling from cells $(i+j-1)\,\mod
n$ and $(i-j+1)\,\mod n$ to $i$ are identical. Since the dynamics of
all cells is given by the same function $f$, this is true for all
$i=1,\ldots,n$. Thus, $F$ is $\kappa$-equivariant.

Suppose now that $F$ is $\kappa$-equivariant. Equality of both sides
of the equivariance condition implies that for all $i=1,\ldots,n$,
the couplings from cells $(i+j-1)\,\mod n$ and $(i-j+1)\,\mod n$ to
$i$ are identical. Since the dynamics of all cells is given by the
same function $f$, the coupling from cell $(i+j-1)\,\mod n$ to cell
$i$ guarantees an identical coupling from cell $i$ to cell
$(i-j+1)\,\mod n$ and this proves (i). But, the coupling from cell
$(i-j+1)\,\mod n$ to $i$ is therefore identical to the coupling from
cell $i$ to cell $(i-j+1)\,\mod n$. Hence there is an identical
two-way coupling between cells $i$ and $(i-j+1)\,\mod n$ which
proves (ii).\qed

\subsubsection{General form of the characteristic equation}
We now focus our attention on delay-coupled cell systems where each
cell is one-dimensional, that is $k=1$. We split the linear and
nonlinear parts of systems~(\ref{eq:Dn-syst-odd1},
\ref{eq:Dn-syst-even1}) and write the result in abstract form:
\[
\dot X=L X_t+H(X_t)
\]
where $X_t\in C([-\tau,0],\R^{n})$, $L:C([-\tau,0],\R^{n})\to
\R^{n}$ is a bounded linear map and $H$ is a nonlinear mapping.
Thus, $L$ is $\D_n$-equivariant, $\eta(\theta)$ is a $n\times n$
$\D_n$-equivariant matrix of bounded variation and
\[
L \phi=\int_{-\tau}^{0} d\eta(\theta)\phi.
\]

\begin{proposition}\label{prop:eta}
The matrix $\eta(\theta)$ is symmetric
($\eta(\theta)=\eta(\theta)^{T}$) with the properties:
\begin{enumerate}
\item for all $j=1,\ldots,n$: $\eta_{jj}(\theta)=p(\theta)$ for some function
$p$,
\item for all $i,k$ with $i\neq k$ then
$\eta_{ki}(\theta)=\eta_{(2+n-k)i}(\theta)=\eta_{k1}(\theta)$.
\end{enumerate}
\end{proposition}

\proof We use Proposition~\ref{prop:Dn-net} to obtain information on
$\eta$. By part (ii), the matrix $\eta(\theta)$ is symmetric. From
the structure of~(\ref{eq:Dn-syst-odd1},\ref{eq:Dn-syst-even1}), we
deduce that for all $j=1,\ldots,n$, $\eta_{jj}(\theta)=p(\theta)$
for some function $p(\theta)$. We denote by $\eta_{ji}(\theta)$ the
element of $\eta$ corresponding to the coupling from cell $j$ to
$i$. Consider $\eta_{j1}(\theta)$ then there is an identical
connection from cell $1$ to cell $2+n-j$ by part (i) and so
$\eta_{j1}(\theta)=\eta_{1(2+n-j)}(\theta)$. By part (ii), the
connection from cell $2+n-j$ to cell $1$ is identical to its
reciprocal and so: $\eta_{j1}(\theta)=\eta_{(2+n-j)1}(\theta)$. By
part (i), we then have
$\eta_{ki}(\theta)=\eta_{(2+n-k)i}(\theta)=\eta_{k1}(\theta)$ since
the connections to cell $i$ are identical to the connections to cell
$1$. \qed

\begin{remark}
This result can be obtained for higher-dimensional cells with a
proof essentially similar to this one, but with a more cumbersome
notation. We decided to restrict ourselves to the one-dimensional
case as this is the one which we study in details in what follows.
\end{remark}

For the purpose of diagonalizing the linear equation, it is
preferable to use $\C^{n}$, the complexification of $\R^{n}$. The
decomposition of $\C^{n}$ by the action of $\D_n=\langle \rho,\kappa
\rangle$ yields an isotypic decomposition of $\C^{n}$, see
Golubitsky {\em et. al.}~\cite{GSS88}. Let $\zeta=e^{2\pi i/n}$ and
define
\[
V_j=\{[v,\zeta^{j}v,\zeta^{2j}v,\ldots,\zeta^{(n-1)j}v]:v\in \R\}
\]
for $j=0,\ldots,n-1$. The subspaces $V_j$ are irreducible
representations of $\D_n$. Recall that all complex irreducible
representations of $\D_n$ are one-dimensional. Then the isotypic
decomposition is
\[
\C^{n}=V_0\oplus V_1\oplus \ldots \oplus V_{n-1}
\]
where $V_0$ is the trivial representation and for $n$ even $V_{n/2}$
is the alternating representation both of which are real
one-dimensional irreducible representations of $\D_n$.

Using Proposition~\ref{prop:eta} and by invariance of the subspace
$V_j$ we have
{\small \[
\eta(\theta)[v,\zeta^{j}v,\zeta^{2j}v,\ldots,\zeta^{(n-1)j}v]=\left(p(\theta)+\sum_{k=2}^{n}
(\zeta^{(k-1)j}+\zeta^{(n+1-k)j})\eta_{k1}(\theta)\right)[v,\zeta^{j}v,\zeta^{2j}v,\ldots,\zeta^{(n-1)j}v].
\]}
Since $\zeta^{n-(k-1)}=\ov{\zeta^{k-1}}$, then for $j=0,\ldots,n-1$
\begin{equation}\label{Aj-theta}
\begin{array}{rcl}
A_j(\theta)&:=&p(\theta)+\displaystyle\sum_{k=2}^{n}
(\zeta^{(k-1)j}+\zeta^{(n+1-k)j})\eta_{k1}(\theta)\\
&=&p(\theta)+\displaystyle\sum_{k=2}^{n}
2\cos(2\pi(k-1)j/n)\eta_{k1}(\theta)\\
\\
&=&\left\{\begin{array}{ll}
    p(\theta)+\displaystyle\sum_{k=2}^{(n+1)/2}
    4\cos(2\pi(k-1)j/n)\eta_{k1}(\theta) & \mbox{$n$ odd}\\
    p(\theta)+2(-1)^{j}\eta_{(1+n/2)1}(\theta)+\displaystyle\sum_{k=2}^{n/2}
    4\cos(2\pi(k-1)j/n)\eta_{k1}(\theta) & \mbox{$n$ even}
    \end{array}\right.
\end{array}
\end{equation}
since $\eta_{k1}=\eta_{(2+n-k)1}$ and
$\cos(2\pi(k-1)j/n)=\cos(2\pi((2+n-k)-1)j/n)$. Note that
$A_j(\theta)=A_{n-j}(\theta)$ for $j=1,\ldots,[n/2]$. The block
diagonalization of $\eta$ is given by the terms $A_{j}(\theta)$ for
$j=0,\ldots,n-1$. Hence, we have
\[
\Delta(\lambda)=\lambda I_n-\int_{-\tau}^{0}
d\eta(\theta)e^{\lambda\theta}=\lambda I_n-
\int_{-\tau}^{0}\diag(dA_{0}(\theta)e^{\lambda\theta},\ldots,dA_{n-1}(\theta)e^{\lambda\theta}).
\]
Let $\Delta_j(\lambda)=\lambda-\int_{-\tau}^{0}
dA_{j}(\theta)e^{\lambda\theta}$, then
\[
\Delta(\lambda)=\diag(\Delta_0(\lambda),\ldots,\Delta_{n-1}(\lambda)).
\]
Therefore, the characteristic equation has, for $n$ even, the
decomposition
\begin{equation}\label{char-eq:Dn-even}
\det\Delta(\lambda)=\det\Delta_{0}(\lambda)\det\Delta_{n/2}(\lambda)\prod_{j=1}^{(n-1)/2}
[\det\Delta_{j}(\lambda)]^{2}=0
\end{equation}
and for $n$ odd
\begin{equation}\label{char-eq:Dn-odd}
\det\Delta(\lambda)=\det\Delta_{0}(\lambda)\prod_{j=1}^{(n-1)/2}
[\det\Delta_{j}(\lambda)]^{2}=0.
\end{equation}

\subsubsection{Application of Main Theorems to the $\D_n$ case}
It is straightforward that Theorem~\ref{thm:main} can be applied to
any of the factors of the characteristic
equations~(\ref{char-eq:Dn-even}) and~(\ref{char-eq:Dn-odd}).

To apply Theorem~\ref{thm:main2} in the $\D_n$ case, we need to
verify that the coefficients $b_{k}^{j}$ in the factors of the
characteristic equation are nonzero and that the nondegeneracy
condition $\det {\cal I}_{B}\neq 0$ is satisfied. In fact, as it is
shown in Section~\ref{sec:main-proof}, $\det {\cal I}_{B}\neq 0$ if
and only if the matrix
\[
{\cal B}:=\left[\begin{array}{cccc} b_1^1 & b_{1+\mu_1}^{1} & \cdots
&
b_{1+\mu_{r-1}}^{1} \\
b_1^2 & b_{1+\mu_1}^2 & \cdots & b_{1+\mu_{r-1}}^2\\
\vdots & \vdots & \ddots & \vdots \\
b_1^r & b_{1+\mu_1}^{r} & \cdots & b_{1+\mu_{r-1}}^{r}
\end{array}\right].
\]
is nonsingular where $\ell_j$ is the number of imaginary eigenvalues
satisfying the $j^{th}$ term of the product of the characteristic
equation~(\ref{char-eq:product}) and $\mu_{j}=\sum_{i=1}^{j} \ell_i$
where $\mu_r=n$ and $\mu_0:=0$. Note that row $j$ of ${\cal B}$
contains coefficients belonging to the $j^{th}$ factor of the
characteristic equation~(\ref{char-eq:product}). The cases $n$ even
and $n$ odd differ significantly and we focus initially on the case
$n$ odd.

\subsection{Critical eigenvalues of $\D_n$ coupled cell
systems}\label{section:critical-Dn} In this section, we apply
Theorem~\ref{thm:main2} to $\D_n$-symmetric coupled cell system
depending on an arbitrary number of finite delays.

\subsubsection{$n$ odd}
The characteristic equation is
\begin{equation}
\det\Delta(\lambda)=\det\Delta_{0}(\lambda)\prod_{j=1}^{(n-1)/2}
[\det\Delta_{j}(\lambda)]^{2}=0.
\end{equation}
We can write
\[
\Delta_{j}(\lambda)=\lambda-F(\lambda)-G_{j}(\lambda)
\]
where
\[
F(\lambda)=\sum_{i=1}^{p} a_i e^{-\lambda\tau_i}
\]
are the terms coming from the internal dynamics of each cell and
\[
G_{j}(\lambda)=\sum_{k=2}^{(n+1)/2}\left[4\cos\left(\dfrac{2\pi(k-1)j}{n}\right)\right]\sum_{t=1}^{m_{k}}
\alpha_{t}^{k} e^{-\lambda s_{t}^{k}}.
\]
are the contributions from the coupling where $m_{k}$ is the number
of delayed terms in the connection from cell $k$ to $1$ and
$\alpha_{t}^{k}$ are the respective coupling coefficients.

\begin{example}
As an example, consider a delay-coupled $\D_5$-symmetric cell. Let
$u_{s}(\theta)=0$ if $\theta=[-\tau,-s]$ and $u_{s}(\theta)=1$ for
$\theta\in (-s,0]$ where $\tau\geq s$ for all delays $s$ and suppose
\[
\eta(\theta)=\left[\begin{array}{ccccc} p(\theta) &
\eta_{21}(\theta) & \eta_{31}(\theta) &
\eta_{41}(\theta) & \eta_{51}(\theta)\\
\eta_{21}(\theta) & p(\theta) & \eta_{51}(\theta) &
\eta_{31}(\theta) & \eta_{41}(\theta)
 \\
\eta_{41}(\theta) & \eta_{21}(\theta) & p(\theta) &
\eta_{51}(\theta) & \eta_{31}(\theta)
  \\
\eta_{31}(\theta) & \eta_{41}(\theta) & \eta_{21}(\theta) &
p(\theta) & \eta_{51}(\theta)
  \\
\eta_{51}(\theta) & \eta_{41}(\theta) & \eta_{31}(\theta) &
\eta_{21}(\theta) & p(\theta)
\end{array}\right]
\]
where
\[
p(\theta)=\dsum_{i=1}^2 a_{i} u_{\tau_i}(\theta),\quad
\eta_{21}(\theta)=\dsum_{\ell=1}^3 \alpha_{\ell}^2
u_{s_{\ell}^2}(\theta) \AND \eta_{31}(\theta)=\dsum_{\ell=1}^2
\alpha_{\ell}^3 u_{s_{\ell}^3}(\theta)
\]
with the conditions $\eta_{41}(\theta)=\eta_{31}(\theta)$ and
$\eta_{51}(\theta)=\eta_{21}(\theta)$ given by
Proposition~\ref{prop:eta} part (2). Then,
\[
F(\lambda)=\sum_{i=1}^2 a_i e^{-\lambda \tau_i}
\]
and
\[
G_{j}(\lambda)=\sum_{k=2}^{3}\left[4\cos\left(\dfrac{2\pi(k-1)j}{n}\right)\right]\sum_{t=1}^{m_{k}}
\alpha_{t}^{k} e^{-\lambda s_{t}^{k}}
\]
where $m_2=3$ and $m_3=2$.
\end{example}

Thus, all coefficients $b_{j}^{k}$ of ${\cal I}_{B}$ are nonzero and
it is convenient to set $b_1^{j}$ to be the coefficient of $a_1$;
that is $b_1^j=1$ for $j=0,1,2,\ldots,(n+1)/2$, and we keep this
convention for the remainder of the paper.

We suppose that the characteristic equation $\Delta(\lambda)=0$ has
purely imaginary roots coming from all factors, then for $r=(n-1)/2$
we have
\[
{\cal B}:=\left[\begin{array}{cccc} b_1^1 & b_{1+\mu_1}^{1} & \cdots
&
b_{1+\mu_{r-1}}^{1} \\
b_1^2 & b_{1+\mu_1}^2 & \cdots & b_{1+\mu_{r-1}}^2\\
\vdots & \vdots & \ddots & \vdots \\
b_1^r & b_{1+\mu_1}^{r} & \cdots & b_{1+\mu_{r-1}}^{r}
\end{array}\right].
\]
and we assign the coefficients $b_{1+\mu_{j}}^{i}$ as follows. We
suppose that the first row corresponds to the factor for the trivial
representation which means that
\[
b_{1+\mu_j}^1=4,\qquad j=1,2,\ldots,(n-1)/2.
\]
Then, we set the remaining coefficients of each row to be equal to
$\left[4\cos\left(\dfrac{2\pi(k-1)j}{n}\right)\right]$ for
$k=2,3,\ldots,(n+1)/2$ where row $j+1$ has the coefficients of
$\Delta_{j}$ for $j=1,2,\ldots,(n-1)/2$. This leads to the matrix
\begin{equation}\label{Dn-B-matrix}
{\cal B}=\left[\begin{array}{cccccc}
 1 & 4 & 4
& \cdots & 4 & 4 \\
1 & 4\cos\left(\dfrac{2\pi}{n}\right) &
4\cos\left(\dfrac{4\pi}{n}\right)
& \cdots & 4\cos\left(\dfrac{(n-3)\pi}{n}\right) & 4\cos\left(\dfrac{(n-1)\pi}{n}\right) \\
1 & 4\cos\left(\dfrac{4\pi}{n}\right) &
4\cos\left(\dfrac{8\pi}{n}\right)
& \cdots & 4\cos\left(\dfrac{2(n-3)\pi}{n}\right) & 4\cos\left(\dfrac{2(n-1)\pi}{n}\right) \\
\\
\vdots & \vdots & \vdots & \vdots & \vdots & \vdots \\
\\
1 & 4\cos\left(\dfrac{(n-1)\pi}{n}\right) &
4\cos\left(\dfrac{2(n-1)\pi}{n}\right) & \cdots &
4\cos\left(\dfrac{(n-3)(n-1)\pi}{2n}\right) &
4\cos\left(\dfrac{(n-1)^2\pi}{2n}\right)
\end{array}\right].
\end{equation}

\vspace{1\baselineskip}

Let $i_1<i_2<\ldots<i_s$ be a set of indices chosen from
$\{0,\ldots,(n-1)/2\}$ defining a combination of factors from the
characteristic equation~(\ref{char-eq:Dn-odd}). We now construct the
$s\times s$ matrix ${\cal B}$ by removing rows and columns
of~(\ref{Dn-B-matrix}) not in the set $\{i_1,i_2,\ldots,i_s\}$.
Suppose that $i_1,\ldots,i_s$ are chosen from $\{1,\ldots,(n-1)/2\}$
then the matrix ${\cal B}$ is symmetric (${\cal B}^{T}={\cal B}$)
and has the form

{\small \begin{equation}\label{matrixB-Dn-a} {\cal
B}=\left[\begin{array}{ccccc} 4\cos\left(\dfrac{2\pi i_1^2
}{n}\right) & 4\cos\left(\dfrac{2\pi i_2 i_1 }{n}\right) & \cdots &
4\cos\left(\dfrac{2\pi i_{s-1} i_1 }{n}\right) &
4\cos\left(\dfrac{2\pi i_s i_1 }{n}\right) \\
\\
\\
4\cos\left(\dfrac{2\pi i_1 i_2 }{n}\right) & 4\cos\left(\dfrac{2\pi
i_2^2 }{n}\right) & \cdots & 4\cos\left(\dfrac{2\pi i_{s-1} i_2
}{n}\right)
& 4\cos\left(\dfrac{2\pi i_s i_2 }{n}\right)\\
\\
\\
\vdots & \vdots & \vdots & \vdots  & \vdots \\
\\
\\
4\cos\left(\dfrac{2\pi i_{1} i_{s-1} }{n}\right) &
4\cos\left(\dfrac{2\pi i_{2} i_{s-1} }{n}\right) & \cdots &
4\cos\left(\dfrac{2\pi i_{s-1}^2 }{n}\right) &
4\cos\left(\dfrac{2\pi i_{s}i_{s-1} }{n}\right)
\\
\\
4\cos\left(\dfrac{2\pi i_1 i_s}{n}\right) & 4\cos\left(\dfrac{2\pi
i_2 i_s}{n}\right) & \cdots & 4\cos\left(\dfrac{2\pi i_{s-1} i_s
}{n}\right) & 4\cos\left(\dfrac{2\pi i_{s}^2}{n}\right).
\end{array}\right].
\end{equation}}

In the other case, $i_1=0$ and the matrix is of the form

{\small \begin{equation}\label{matrixB-Dn-b} {\cal
B}=\left[\begin{array}{ccccc} 1 & 4 & \cdots & 4 &
4 \\
\\
\\
1 & 4\cos\left(\dfrac{2\pi i_2^2}{n}\right) & \cdots &
4\cos\left(\dfrac{2\pi i_{s-1} i_2-1 }{n}\right)
& 4\cos\left(\dfrac{2\pi i_s i_2 }{n}\right)\\
\\
\\
\vdots & \vdots & \vdots & \vdots  & \vdots \\
\\
\\
1 & 4\cos\left(\dfrac{2\pi i_{2} i_{s-1} }{n}\right) & \cdots &
4\cos\left(\dfrac{2\pi i_{s-1} i_{s-1} }{n}\right) &
4\cos\left(\dfrac{2\pi i_{s} i_{s-1} }{n}\right)
\\
\\
1 & 4\cos\left(\dfrac{2\pi i_2 i_s}{n}\right) & \cdots &
4\cos\left(\dfrac{2\pi i_{s-1}^2 }{n}\right) &
4\cos\left(\dfrac{2\pi i_{s}^2}{n}\right).
\end{array}\right].
\end{equation}}
We can now state our result.

\begin{theorem}\label{Thm:Dn-odd}
Consider a linear $\D_n$-symmetric coupled cell system with $n$ odd
depending on $k$ delays $\tau_1,\ldots,\tau_k$ and let
$i_1<i_2<\ldots<i_s$ be indices chosen from $\{0,\ldots,(n-1)/2\}$
defining a combination of factors from the characteristic
equation~(\ref{char-eq:Dn-odd}). We assume that the matrix ${\cal
B}$ given by~(\ref{matrixB-Dn-a}) or~(\ref{matrixB-Dn-b}) is
nonsingular. Suppose
\[
\omega_1^1,\ldots,\omega_{\ell_{i_1}}^{1},\omega_1^2,\ldots,\omega_{\ell_{i_2}}^{2}\ldots,\omega_{1}^{s},\ldots,\omega_{\ell_{i_s}}^{s}
\]
are positive and linearly independent over the rationals, where
$\ell_{i_1}+\cdots+\ell_{i_s}=k$. Then there exists
$\tau_1>0,\ldots,\tau_k>0$ and real coefficients $a_i$ such that for
all $m=1,\ldots,s$
\[
\Delta_{i_m}(\lambda)=0
\]
has roots $i\omega_{\ell}^{m}$ for $\ell=1,\ldots,\ell_{i_m}$.
\end{theorem}

\proof Since $n$ is odd, the coefficients
$$b_k^{j}=4\cos\left(\frac{2\pi(k-1)j}{n}\right)$$ are nonzero for all
$k=2,\ldots,(n+1)/2$ and $j=0,\ldots,(n-1)/2$. Because ${\cal B}$ is
assumed nonsingular, Theorem~\ref{thm:main2} applies and the result
is obtained. \qed.

The condition that ${\cal B}$ is nonsingular does not always hold as
we show in the case $s=2$. Consider the matrix~(\ref{matrixB-Dn-b})
with $n=9$ so that $i_2\in \{1,2,3,4\}$. Choosing $i_2=3$ we have
the singular matrix
\[
{\cal B}=\left(\begin{array}{cc} 1 &
4 \\
\\ 1 &
4
\end{array}\right).
\]

We now show that ${\cal B}$ is nonsingular in the case $s=2$ if the
matrix is given by~(\ref{matrixB-Dn-a}); that is,
\[
{\cal B}=\left(\begin{array}{cc} 4\cos\left(\dfrac{2\pi i_1^2
}{n}\right) &
4\cos\left(\dfrac{2\pi i_1 i_2}{n}\right) \\
\\ 4\cos\left(\dfrac{2\pi i_2 i_1 }{n}\right) &
4\cos\left(\dfrac{2\pi i_2^2 }{n}\right)
\end{array}\right).
\]
We compute
\[
\begin{array}{rcl}
\det {\cal B}&=&16 \left[ \cos\left(\dfrac{2\pi i_1^2
}{n}\right)\cos\left(\dfrac{2\pi i_2^2
}{n}\right)-\cos\left(\dfrac{2\pi i_1 i_2}{n}\right)^2\right]\\
\\
&=&8\left[\cos\left(\dfrac{2\pi
(i_1^2+i_2^2)}{n}\right)+\cos\left(\dfrac{2\pi
(i_1^2-i_2^2)}{n}\right)-\cos\left(\dfrac{4\pi i_1 i_2}{n}\right)-1
\right].
\end{array}
\]
We do a few cases explicitly. First, the case $n=3$ is not relevant
since $i_1<i_2$, $(n-1)/2=1$ and $i_1\neq 0$. We do the case $n=5$
where we must have $i_1=1$ and $i_2=2$. This implies that
$i_1^2+i_2^2=5$ and so
\[
\det{\cal
B}=8\left[\cos\left(\dfrac{4\pi}{5}\right)-\cos\left(\dfrac{8\pi}{5}\right)\right]\neq
0.
\]
We now turn to the general case and show that the determinant cannot
vanish. Because the three cosines are projections of $n^{th}$ roots
of unity on the real axis for $n$ odd, then
\[
\cos\left(\dfrac{2\pi (i_1^2+i_2^2)}{n}\right)+\cos\left(\dfrac{2\pi
(i_1^2-i_2^2)}{n}\right)-\cos\left(\dfrac{4\pi i_1
i_2}{n}\right)\neq 1.
\]
So if the determinant is to vanish, one of the cosines must be equal
to $1$. Since $i_1<i_2$ there is only one option and we must have
$i_1^2+i_2^2=n$. Thus, $i_1^2=n-i_2^2$ and
\[
\cos\left(\dfrac{2\pi (i_1^2-i_2^2)}{n}\right)=\cos\left(\dfrac{2\pi
(n-2i_2^2)}{n}\right)=\cos\left(\dfrac{4\pi i_2^2}{n}\right).
\]
If
\[
\cos\left(\dfrac{4\pi i_2^2}{n}\right)-\cos\left(\dfrac{4\pi i_1
i_2}{n}\right)=0,
\]
this would impliy $i_1=i_2$, but we know that $i_1<i_2$ and so
$\det{\cal B}$ cannot vanish. We summarize this result in the next
theorem.

\begin{theorem}
Consider a linear $\D_n$-symmetric coupled cell system with $n$ odd
depending on $k$ delays $\tau_1,\ldots,\tau_k$ and let $i_1<i_2$ be
indices chosen from $\{1,\ldots,(n-1)/2\}$ defining a combination of
factors from the characteristic equation~(\ref{char-eq:Dn-odd}).
Suppose
\[
\omega_1^1,\ldots,\omega_{\ell_{i_1}}^{1},\omega_1^2,\ldots,\omega_{\ell_{i_2}}^{2}
\]
are positive and linearly independent over the rationals, where
$\ell_{i_1}+\ell_{i_2}=k$. Then there exists
$\tau_1>0,\ldots,\tau_k>0$ and real coefficients
$a_1\in\R,\ldots,a_p\in \R$ such that for $m=1$ and $m=2$,
\[
\Delta_{i_m}(\lambda)=0
\]
has roots $i\omega_{\ell}^{m}$ for $\ell=1,\ldots,\ell_{i_m}$.
\end{theorem}

\subsection{$n$ even} For $\D_n$-symmetric systems with $n$ even,
Theorem~\ref{thm:main2} cannot be applied directly because the
condition $b_{k}^{j}\neq 0$ for all indices $j,k$ is not always
satisfied. For instance, if $n=4$ then
\[
A_{j}(\theta)=p(\theta)+2(-1)^{j}\eta_{31}(\theta)+4\cos(2\pi
j/4)\eta_{21}(\theta)
\]
and we have
\[
A_1(\theta)=A_3(\theta)=p(\theta)-2\eta_{31}(\theta).
\]
Hence, $b_{2}^{1}=b_{2}^{3}=0$.

\section{Proof of Theorem~\ref{thm:main2}}\label{sec:main-proof}
Before we present the proof of Theorem~\ref{thm:main2}, we describe
in the next lemma the form of the matrix ${\cal I}_{B}$ which
appears in the proof and compute its determinant.
\begin{lemma}
Let $\ell_1,\ldots,\ell_r$ be positive integers and define
$\mu_{j}=\sum_{i=1}^{j} \ell_i$ where $\mu_r=n$ and $\mu_0:=0$.
Consider the $n\times n$ matrix
\[
{\cal I}_{B}:=\left[ A_1 \cdots A_j \cdots A_r \right]^{T}
\]
where
\[
A_j=\left[
\begin{array}{cccccccccc}
b_{1}^{j} & \cdots & b_{\mu_{j-1}}^{j} & b_{1+\mu_{j-1}}^{j} &  b_{2+\mu_{j-1}}^{j} & \cdots &  b_{\mu_j}^{j} & b_{\mu_{j}+1}^{j} & \cdots & b_{n}^{j} \\
b_{1}^{j} & \cdots & b_{\mu_{j-1}}^{j} & b_{1+\mu_{j-1}}^{j} &  b_{2+\mu_{j-1}}^{j} & \cdots & - b_{\mu_j}^{j} & b_{\mu_{j}+1}^{j} & \cdots & b_{n}^{j}\\
\vdots & \cdots &\vdots & \vdots &\vdots & \cdots & \vdots & \vdots &\cdots & \vdots \\
b_{1}^{j} & \cdots & b_{\mu_{j-1}}^{j}& b_{1+\mu_{j-1}}^{j} &
-b_{2+\mu_{j-1}}^{j} & \cdots & -b_{\mu_j}^{j} & b_{\mu_j+1}^{j} &
\cdots & b_{n}^{j}
\end{array}\right]
\]
is a $\ell_j\times n$ matrix and all elements are nonzero. Then,
\[
\det {\cal I}_{B}=\pm \prod_{j=1}^{r} \left[(-2)^{\ell_j-1}
\prod_{s=2}^{\ell_j} b_{s+\mu_{j-1}}^{j} \right]\det {\cal B}
\]
where
\[
{\cal B}:=\left[\begin{array}{cccc} b_1^1 & b_{1+\mu_1}^{1} & \cdots
&
b_{1+\mu_{r-1}}^{1} \\
b_1^2 & b_{1+\mu_1}^2 & \cdots & b_{1+\mu_{r-1}}^2\\
\vdots & \vdots & \ddots & \vdots \\
b_1^r & b_{1+\mu_1}^{r} & \cdots & b_{1+\mu_{r-1}}^{r}
\end{array}\right].
\]
\end{lemma}

\proof Substitute row $k$, denoted by $R_k$, of
\[
A_j=\left[
\begin{array}{cccccccccc}
b_{1}^{j} & \cdots & b_{\mu_{j-1}}^{j} & b_{1+\mu_{j-1}}^{j} &  b_{2+\mu_{j-1}}^{j} & \cdots &  b_{\mu_j}^{j} & b_{\mu_{j}+1}^{j} & \cdots & b_{n}^{j} \\
b_{1}^{j} & \cdots & b_{\mu_{j-1}}^{j} & b_{1+\mu_{j-1}}^{j} &  b_{2+\mu_{j-1}}^{j} & \cdots & - b_{\mu_j}^{j} & b_{\mu_{j}+1}^{j} & \cdots & b_{n}^{j}\\
\vdots & \cdots &\vdots & \vdots &\vdots & \cdots & \vdots & \vdots &\cdots & \vdots \\
b_{1}^{j} & \cdots & b_{\mu_{j-1}}^{j}& b_{1+\mu_{j-1}}^{j} &
-b_{2+\mu_{j-1}}^{j} & \cdots & -b_{\mu_j}^{j} & b_{\mu_j+1}^{j} &
\cdots & b_{n}^{j}
\end{array}\right]
\]
for $k=2,\ldots,\ell_j$ by $R_k-R_1$. The matrix $A_j$ becomes
{\small \[ \tilde{A}_j:=\left[
\begin{array}{cccccccccccc}
b_{1}^{j} & \cdots & b_{\mu_{j-1}}^{j} & b_{1+\mu_{j-1}}^{j} &
b_{2+\mu_{j-1}}^{j}
& b_{3+\mu_{j-1}}^j & \cdots & b_{-1+\mu_{j}}^j & b_{\mu_j}^{j} & b_{\mu_{j}+1}^{j} & \cdots & b_{n}^{j} \\
0 & \cdots & 0 & 0 &  0 & 0 & \cdots & 0 & -2 b_{\mu_j}^{j} & 0 & \cdots & 0\\
0 & \cdots & 0 & 0 &  0 & 0 & \cdots & -2 b_{-1+\mu_j}^j & -2
b_{\mu_j}^j & 0 & \cdots & 0\\
\vdots & \cdots &\vdots & \vdots &\vdots & \vdots & \cdots & \vdots
& \vdots & \vdots
&\cdots & \vdots \\
0 & \cdots & 0 & 0 & 0 & -2b_{3+\mu_{j-1}}^j & \cdots &
-2b_{-1+\mu_{j}}^j & -2
b_{\mu_j}^j & 0 & \cdots & 0\\
 0 & \cdots & 0 & 0 & -2b_{2+\mu_{j-1}}^{j} & -2 b_{3+\mu_{j-1}}^j & \cdots & -2
b_{-1+\mu_j}^j & -2b_{\mu_j}^{j} & 0 & \cdots & 0
\end{array}\right]
\]}
We compute the determinant of ${\cal I}_B$ by cofactor expansion
starting with row $2$ of $\tilde{A}_j$ which contains a unique
nonzero element $-2 b_{\mu_j}^j$. Denote by $C_{ij}$ the
$(i,j)$-cofactor matrix. The row $2+\mu_{j-1}$ of
$C_{(2+\ell_{j-1},\mu_j)}$ has a unique nonzero element
$-2b_{-1+\mu_{j}}$ and we perform a cofactor expansion along this
row. The row $2+\mu_{j-1}$ of this new cofactor matrix also has a
unique nonzero element $-2 b_{-2+\mu_{j}}^j$ and we proceed with the
same process removing successively columns $3+\mu_{j-1}$ to $\mu_j$
(and the appropriate rows) until the cofactor matrix has only two
rows corresponding to the original $\tilde{A}_j$ matrix and the
second row has the unique nonzero element $-2 b_{2+\mu_{j-1}}$ which
is used to perform a cofactor expansion. Performing this process
successively on each matrix $\tilde{A}_j$ for $j=1,\ldots,r$, leaves
as a cofactor matrix the $r\times r$ matrix ${\cal B}$ defined in
the statement. The formula in the lemma is written using
$\mu_j=\mu_{j-1}+\ell_j$ and so the lemma is proved.\qed

\noindent We are now ready to prove our Theorem~\ref{thm:main2}.

\vspace{1\baselineskip}

\proofof{Theorem~\ref{thm:main2}} A necessary and sufficient
condition for the conclusion of the theorem to hold is that the
following algebraic system of $2n$ equations has a solution in the
$2n$ unknowns $(\tau_1,\tau_2,\ldots,\tau_n,a_1,a_2,\ldots,a_n)$:
\begin{equation}
\left\{\begin{array}{ccc}
{\displaystyle\sum_{k=1}^n\,a_k\,(b_k^{1} e^{-i\omega_{\ell}^{1}\tau_k})}&=&i\omega_{\ell}^{1},\,\,\,\,\,\,\ell=1,\ldots,\ell_1\\[0.15in]
{\displaystyle\sum_{k=1}^n\,a_k\,(b_k^{2} e^{-i\omega_{\ell}^{1}\tau_k})}&=&i\omega_{\ell}^{2},\,\,\,\,\,\,\ell=1,\ldots,\ell_2\\[0.15in]
\vdots \\
{\displaystyle\sum_{k=1}^n\,a_k\,(b_k^{r}
e^{-i\omega_{\ell}^{1}\tau_k})}&=&i\omega_{\ell}^{r},\,\,\,\,\,\,\ell=1,\ldots,\ell_r
\end{array}\right.
\end{equation}
\begin{equation}
\left\{\begin{array}{ccc} {\displaystyle\sum_{k=1}^n\,a_k\,(b_k^{1}
e^{i\omega_{\ell}^{1}\tau_k})}&=&-i\omega_{\ell}^{1},\,\,\,\,\,\,k=1,\ldots,\ell_1\\[0.15in]
{\displaystyle\sum_{k=1}^n\,a_k\,(b_k^{2} e^{i\omega_{\ell}^{1}\tau_k})}&=&-i\omega_{\ell}^{2},\,\,\,\,\,\,\ell=1,\ldots,\ell_2\\[0.15in]
\vdots\\
{\displaystyle\sum_{k=1}^n\,a_k\,(b_k^{r}
e^{i\omega_{\ell}^{1}\tau_k})}&=&-i\omega_{\ell}^{r},\,\,\,\,\,\,\ell=1,\ldots,\ell_r
\end{array}\right.
\label{chareq-1}
\end{equation}
Although (\ref{chareq-1}) is in complex form, system
(\ref{chareq-1}) is equivalent to a system of $2n$ {\it real}
equations. This fact is taken for granted throughout the sequel,
even though we continue to use complex notation.

Let
\[
\omega=(\omega_{1}^{1},\ldots,\omega_{\ell_1}^{1},\ldots,\omega_{1}^{r},\ldots,\omega_{\ell_r}^{r}).
\]
It is useful to use the following matrix notation for
(\ref{chareq-1})
\begin{equation}
\left(\begin{array}{c}P(\tau;\omega)\\[0.1in]P(-\tau;\omega)\end{array}\right)A^T=\left(\,\begin{array}{c}i\omega^T\\[0.1in]-i\omega^T\end{array}\,\right)
\label{chareq_matrix-1}
\end{equation}
where $A=(a_1,\ldots,a_n)$, superscript $T$ denotes transpose, and
$P(\tau;\omega)=P(\tau_1,\ldots,\tau_n;\omega)$ is the $n\times n$
matrix of the form
\[
P(\tau;\omega)=\left[\begin{array}{c} P_1(\tau;\omega) \\
P_2(\tau;\omega)\\ \vdots \\ P_r(\tau;\omega)\end{array}\right]
\]
whose entry at block $j$, row $\ell$ and column $k$ is
\[
\left[\,P_{j}(\tau;\omega)\,\right]_{\ell k}=b_k^{j}
e^{-i\omega_{\ell}^{j}\tau_k}.
\]
Note that $\overline{P(\tau;\omega)}=P(-\tau;\omega)$.

\vspace*{0.25in} Recall that $\mathrm{T}=(\mathbb{S}^1)^n$ and
$\mathrm{V}=\mathrm{T}^n$. Consider the following mapping associated
to (\ref{chareq_matrix-1}):
\[
F:\mathrm{V}\times\mathbb{R}^n\longmapsto\mathbb{R}^{2n}
\]
defined as follows.
\begin{equation}
F(\Phi^1,\ldots,\Phi^{r},A;\omega)=\left(\begin{array}{c}\widetilde{P}(\Phi^1,\ldots,\Phi^r)\\[0.1in]
\widetilde{P}(-\Phi^1,\ldots,-\Phi^r)\end{array}\right)A^T-i\left(\begin{array}{c}\omega^T\\[0.1in]-\omega^T\end{array}\right),
\label{real_suspended-1}
\end{equation}
where $A$ and $\omega$ are as previously defined and
\[
\widetilde{P}(\Phi^1,\ldots,\Phi^r)=\left[\begin{array}{c} \widetilde{P}_1(\Phi^1) \\
\widetilde{P}_2(\Phi^2)\\ \vdots \\
\widetilde{P}_r(\Phi^r)\end{array}\right]
\]
with
\[
\left[\,\widetilde{P}_{j}(\Phi^1,\ldots,\Phi^n)\,\right]_{\ell
k}=b_{k}^{j} e^{-i\varphi^j_{\ell k}}.
\]
for $j=1,\ldots,r$, $\ell=1,\ldots,\ell_j$ and $k=1,\ldots,n$. The
definition of $\widetilde{P}$ uses the following coordinates of
$\mathrm{V}$. Let $\mathrm{T}_j=(\mathbb{S}^1)^{\ell_j}$ for
$j=1,\ldots,r$, then
\[
(\Phi^1,\ldots,\Phi^r)\in\mathrm{T}_{j}^n\times\cdots\times
\mathrm{T}_{j}^n
\]
where
\[
\Phi^{j}=(\Phi_1^j,\ldots,\Phi_n^j) \AND
\Phi_{k}^{j}=(\varphi_{1k}^{j},\ldots,\varphi_{\ell_j
k}^{j})^{T}\in\mathrm{T}_j.
\]
Let $\Psi_j=\Phi_{n}^{j}$ and $\Psi=(\Psi_1,\ldots,\Psi_r)$. We use
the notation
\[
\Phi=(\Phi_{o}^1,\ldots,\Phi_{o}^{n-1})
\]
where
\[
\Phi_{o}^{j}=(\Phi_1^j,\ldots,\Phi_{n-1}^j)
\]
so that the mapping $F$ in (\ref{real_suspended-1}) is written as
$F(\Phi,\Psi,A;\omega)$.

\vspace*{0.25in} \noindent {\bf Definition}: If
$\{e_1,\ldots,e_{\ell}\}$ denotes the canonical basis of vectors in
$\mathbb{R}^{\ell}$, we define the vectors $v_1,\ldots,v_{\ell}$ by
${\displaystyle v_1=e_1+\cdots+e_{\ell}}$, and for
$k=2,\ldots,{\ell}$,
\[
v_k=v_1-\sum_{m=0}^{\ell-2}\,2\,e_{\ell-m}.
\]
By construction, the set $\{v_1,\ldots,v_{\ell}\}$ is linearly
independent and so the $\ell\times \ell$ matrix ${\cal I}^{\ell}$,
whose $j^{\mbox{\small th}}$ column is the vector $v_j^T$, is
invertible.

We also define the $\ell\times \ell$ invertible matrix ${\cal U}_j$
to be the diagonal matrix whose $k^{\mbox{\small th}}$ diagonal
element is the $k^{\mbox{\small th}}$ component of the vector $v_j$
(in particular, ${\cal U}_1$ is the identity matrix). Note also that
${\cal U}_j^{-1}={\cal U}_j$, $j=1,\ldots,r$.

We define $\mu_0:=0$, $\mu_{j}:=\sum_{i=1}^{j} \ell_i$,
\[
\Theta_{\ell_j}:=(\Phi_{1+\mu_{j-1}}^{j},\ldots,\Phi_{\mu_j}^{j})
\]
We use the following base point in
$\mathrm{V}=\mathrm{V}_{\Phi}\times \mathrm{V}_{\Psi}$. For
$j=1,\ldots,r$ define $\hat{\Phi}_j$ be the point given by
\[
\hat{\Theta}_{\ell_j}=\frac{-\pi}{2}((v_1,\ldots,v_{\ell_{j}-1}),v_{\ell_j})
\]
and $\hat{\Phi}_{i}^{j}=-\frac{\pi}{2}v_1$ for
$i\not\in\{\mu_{1+\mu_{j-1}},\ldots,\mu_{j}\}$. In particular,
$\hat{\Phi}_{\mu_{j}}^{j}=-\frac{\pi}{2}v_{\ell_j}$.

We now evaluate $\widetilde{P}(\widehat{\Phi},\widehat{\Psi})$ by
computing $\widetilde{P}_j(\widehat{\Phi},\widehat{\Psi})$ for
$j=1,\ldots,r$:
\[
\widetilde{P}_j(\widehat{\Phi},\widehat{\Psi})=i\left[
\begin{array}{cccccccccc}
b_{1}^{j} & \cdots & b_{\mu_{j-1}}^{j} & b_{1+\mu_{j-1}}^{j} &  b_{2+\mu_{j-1}}^{j} & \cdots &  b_{\mu_j}^{j} & b_{\mu_{j}+1}^{j} & \cdots & b_{n}^{j} \\
b_{1}^{j} & \cdots & b_{\mu_{j-1}}^{j} & b_{1+\mu_{j-1}}^{j} &  b_{2+\mu_{j-1}}^{j} & \cdots & - b_{\mu_j}^{j} & b_{\mu_{j}+1}^{j} & \cdots & b_{n}^{j}\\
\vdots & \cdots &\vdots & \vdots &\vdots & \cdots & \vdots & \vdots &\cdots & \vdots \\
b_{1}^{j} & \cdots & b_{\mu_{j-1}}^{j}& b_{1+\mu_{j-1}}^{j} &
-b_{2+\mu_{j-1}}^{j} & \cdots & -b_{\mu_j}^{j} & b_{\mu_j+1}^{j} &
\cdots & b_{n}^{j}
\end{array}\right].
\]
Thus,
\[
\widetilde{P}(\widehat{\Phi},\widehat{\Psi})=\left[\begin{array}{c}
\widetilde{P}_1(\widehat{\Phi},\widehat{\Psi})\\[0.15in]
\vdots\\
\widetilde{P}_j(\widehat{\Phi},\widehat{\Psi})\\[0.15in]
\vdots\\
\widetilde{P}_r(\widehat{\Phi},\widehat{\Psi})
\end{array}\right]:=i{\cal I}_{B}
\]
where ${\cal I}_{B}$ is invertible by assumption. In particular,
$\widetilde{P}(-\widehat{\Phi},-\widehat{\Psi})=-i{\cal I}_{B}$. We
define
\[
\widehat{A}^T\equiv (\hat{a}_1,\ldots,\hat{a}_n)^T={\cal
I}_{B}^{-1}\omega^T
\]
which leads to the solution:
\[
F(\widehat{\Phi},\widehat{\Psi},\widehat{A};\omega)=0.
\]
Because the $\omega_j$ are rationally independent, it follows that
the components $\hat{a}_k$ of $\widehat{A}$ are all non-zero.
Furthermore, it is easy to compute that
\[
J\equiv
D_{(\Psi,A)}F(\widehat{\Phi},\widehat{\Psi},\widehat{A};\omega)=\left(\begin{array}{cc}
\hat{U} & i {\cal I}_{B} \\ \hat{U} & -i {\cal I}_B
\end{array}\right)
\]
where
\[
\hat{U}=\diag(\hat{a}_{n} b_{n}^1 {\cal U}_{1}^1,\ldots,\hat{a}_{n}
b_{n}^{r-1} {\cal U}_{1}^{r-1},\hat{a}_{n} b_{n}^r {\cal
U}_{\ell_r}^r)
\]
is a $n\times n$ matrix with diagonal blocks of dimensions
$\ell_1\times \ell_1$ to $\ell_r\times \ell_r$. We compute also
\[
K\equiv D_{\Phi}F(\widehat{\Phi},\widehat{\Psi},\widehat{A};\omega)=
\left(\,\begin{array}{c} \hat{K} \\ \hat{K}\end{array}\,\right)
\]
where
\[
\hat{K}=\diag(\hat{K}_1,\ldots,\hat{K}_r)
\]
is a $n\times (n-1)n$ matrix where the block $\hat{K}_j$ has
dimensions $\ell_j\times (n-1)\ell_j$ and is of the form
\[
\begin{array}{l}
\hat{K}_j=\left(\begin{array}{cccccc}\hat{a}_1\,b_1^j\,{\cal U}_1^j
& \cdots & \,\hat{a}_{1+\mu_{j-1}}\,b_{1+\mu_{j-1}}^j\,{\cal U}_1^j
& \,\hat{a}_{2+\mu_{j-1}}\,b_{2+\mu_{j-1}}^j\,{\cal U}_2^j &\cdots &
\,\hat{a}_{\mu_j}\,b_{\mu_j}^j\,{\cal U}_{\ell_j}^j
\end{array}\right. \\
\left.\begin{array}{ccc} \hat{a}_{1+\mu_{j}} b_{1+\mu_{j}}^j {\cal
U}_1^j & \cdots & \hat{a}_{n-1} b_{n-1}^j {\cal U}_1^j
\end{array}\right).
\end{array}
\]

The matrix $J$ is invertible and its inverse is
\[
J^{-1}=\left(\,\begin{array}{cc} \frac{1}{2}\hat{U}^{-1} &
\frac{1}{2}\hat{U}^{-1} \\[0.1in] \frac{-i}{2}{\cal I}_{B}^{-1} &
\frac{i}{2}{\cal I}_{B}^{-1}
\end{array}\,\right).
\]
By the implicit function theorem, there exists a neighborhood ${N}$
of $\widehat{\Phi}$ in $\mathrm{V}_{\Phi}$ and a unique smooth
function
\[
\begin{array}{lll}
G:{N}&\longmapsto &\mathrm{T}\times\mathbb{R}^n\\[0.15in]
G: \Phi&\longmapsto&G(\Phi)=(G_{\Psi}(\Phi),G_A(\Phi))
\end{array}
\]
such that
\[
G(\widehat{\Phi})=(\widehat{\Psi},\widehat{A})
\]
and
\begin{equation}
F(\Phi,G(\Phi);\omega)\equiv 0,\,\,\,\forall\,\Phi\in {N}.
\label{ift-1}
\end{equation}

Implicit differentiation of (\ref{ift-1}) yields that
\begin{equation}
\begin{array}{lll}
DG(\widehat{\Phi})&=&\left(\,\begin{array}{c}DG_{\Psi}(\widehat{\Phi})\\[0.1in]
DG_A(\widehat{\Phi})\end{array}\,\right)=
-J^{-1}K\\[0.25in]
&=&\left(\,\begin{array}{c} \diag({\cal M}_1,\ldots,{\cal M}_r)\\[0.1in]
* \end{array}\,\right),
\end{array}
\end{equation}
where $*$ is not important for our purposes and the first component
is a $n\times (n-1)n$ matrix composed of $n-1$ block matrices ${\cal
M}_{j}=$
\[
\begin{array}{l}
\left(\,\begin{array}{cccccccc} -\frac{\hat{a}_{1}
b_{1}^j}{\hat{a}_{\mu_j} b_{\mu_j}^j}\,({\cal U}_1^{j})^2 & \cdots &
-\frac{\hat{a}_{1+\mu_{j-1}} b_{1+\mu_{j-1}}^j}{\hat{a}_{\mu_j}
b_{\mu_j}^j}\,({\cal U}_1^{j})^2 & -\frac{\hat{a}_{2+\mu_{j-1}}
b_{2+\mu_{j-1}}^j}{\hat{a}_{\mu_{j}} b_{\mu_{j}}^j}\,{\cal
U}_1^{j}{\cal U}_2^j &\cdots
 \end{array} \right.\\
\left.\begin{array}{cccc} -\frac{\hat{a}_{\mu_{j}} b_{\mu_{j}}^j
}{\hat{a}_{\mu_{j}} b_{\mu_{j}}^j} \,{\cal U}_1^{j} {\cal
U}_{\ell_j}^{j} & -\frac{\hat{a}_{1+\mu_{j}}
b_{1+\mu_{j}}^j}{\hat{a}_{\mu_j} b_{\mu_j}^j}\,({\cal U}_1^{j})^2 &
\cdots & -\frac{\hat{a}_{n-1} b_{n-1}^j}{\hat{a}_{\mu_j}
b_{\mu_j}^j}\,({\cal U}_1^{j})^2
\end{array}\,\right)
\end{array}
\]
of dimension $\ell_j\times (n-1)\ell_j$ where $j=1,\ldots,r-1$.
Recall that $\mu_r=n$, so that we have
\[
\begin{array}{c}
{\cal M}_r=\left(\,\begin{array}{ccccc} -\frac{\hat{a}_{1}
b_{1}^r}{\hat{a}_{\mu_r} b_{\mu_r}^r}\,{\cal U}_{\ell_r}^{r}\,{\cal
U}_1^{r} & \cdots & -\frac{\hat{a}_{\mu_{r-1}}
b_{\mu_{r-1}}^r}{\hat{a}_{\mu_r} b_{\mu_r}^r}\,{\cal
U}_{\ell_r}^{r}\,{\cal U}_1^{r} & -\frac{\hat{a}_{1+\mu_{r-1}}
b_{1+\mu_{r-1}}^r}{\hat{a}_{\mu_r} b_{\mu_r}^r}\,{\cal
U}_{\ell_r}^{r}\,{\cal U}_1^{r}  \end{array} \right.\\
\left.\begin{array}{ccc} -\frac{\hat{a}_{2+\mu_{r-1}}
b_{2+\mu_{r-1}}^r}{\hat{a}_{\mu_{r}} b_{\mu_{r}}^r}\,{\cal
U}_{\ell_r}^r\,{\cal U}_2^r&\cdots & -\frac{\hat{a}_{n-1} b_{n-1}^r
}{\hat{a}_{\mu_{r}} b_{\mu_{r}}^r} \,{\cal U}_{\ell_r}^{r}\,{\cal
U}_{\ell_{r}-1}^{r}
\end{array}\,\right)
\end{array}
\]
Consequently,
\begin{equation}
DG_{\Psi}(\widehat{\Phi})=\diag({\cal M}_1,\ldots,{\cal M}_r)
\end{equation}
is nonsingular and it follows that the mapping
\[
G_{\Psi}:N\longrightarrow \mathrm{T}
\]
is regular at $\widehat{\Phi}$.

Consider the following $n-1$ vectors in $
(\R^{\ell_1})^{n-1}\times(\R^{\ell_2})^{n-1}\times \cdots\times
(\R^{\ell_r})^{n-1}\simeq(\,\mathbb{R}^n\,)^{n-1}: $
\[
\begin{array}{lll}
W_1&=&(\omega^1,0,\ldots,0;\omega^2,0,\ldots,0;\ldots;\omega^r,0,\ldots,0)\\[0.1in]
W_2&=&(0,\omega^1,\ldots,0;0,\omega^2,\ldots,0;\ldots,0,\omega^r,\ldots,0)\\[0.1in]
&\vdots&\\[0.1in]
W_{n-2}&=&(0,\ldots,\omega^1,0;0,\ldots,\omega^2,0;\ldots;0,\ldots,\omega^r,0)\\[0.1in]
W_{n-1}&=&(0,\ldots,0,\omega^1;0,\ldots,0,\omega^2;\ldots,0,\ldots,0,\omega^r),
\end{array}
\]
where $0$ represents the $0$ vector in the respective space
$\mathbb{R}^{\ell_j}$, and we recall that
\[
(\omega^{1},\ldots,\omega^{r})=(\omega_{1}^{1},\ldots,\omega_{\ell_1}^{1},\ldots,\omega_{1}^{r},\ldots,\omega_{\ell_r}^{r}).
\]
The set $\{W_1,\ldots,W_{n-1}\}$ is linearly independent, so for any
$\Phi\in\,\mathrm{V}_{\Phi}$, the set
\[
{\cal S}_{\Phi}=\{\,\Phi+\sum_{j=1}^{n-1}\sigma_j\,W_j\,\,(\mbox{\rm
  mod}\,\,\mathrm{V}_{\Phi})\,\,|\,\,\,0\leq |\sigma_j| < <
  1,\,\,j=1,\ldots,n\,\}
\]
is a small $n-1$-dimensional surface through $\Phi$ in
$\mathrm{V}_{\Phi}$. We are interested in showing that for $\Phi$
close enough to $\widehat{\Phi}$ in $N$, the image of ${\cal
S}_{\Phi}$ by $G_{\Psi}$ in $\mathrm{V}_{\Psi}$ is transverse to the
integral curves of the vector field $\dot{\Psi}=\omega$. To show
this, we consider the function
\[
{\cal T}: N\longrightarrow \mathbb{R}
\]
defined by
\begin{equation}
{\cal T}(\Phi)=\mbox{\rm det}\left(\,DG_{\Psi}(\Phi)\cdot
  W_1^T\,\,\,\,\,\,\,DG_{\Psi}(\Phi)\cdot
  W_2^T\,\,\,\,\,\,\,\cdots\,\,\,\,\,\,\,DG_{\Psi}(\Phi)\cdot
  W_{n-1}^T\,\,\,\,\,\,\,\omega^T\,\right)
\end{equation}
and recalling that $({\cal U}_{j}^{i})^{2}=I$ for all $j,i$ we
compute
\begin{equation}
\begin{array}{rcl}
{\cal T}(\hat{\Phi})&=&\mbox{\rm
det}\left(\,DG_{\Psi}(\hat{\Phi})\cdot
W_1^T\,\,\,\,\,\,\,DG_{\Psi}(\hat{\Phi})\cdot
W_2^T\,\,\,\,\,\,\,\cdots\,\,\,\,\,\,\,DG_{\Psi}(\hat{\Phi})\cdot
W_{n-1}^T\,\,\,\,\,\,\,\omega^T\,\right)\\
\\
&=&\det\left(\alpha_{jk}\right)
\end{array}
\end{equation}
where $j=1,\ldots,r$, $k=1,\ldots,n$. The elements of the matrix
$(\alpha_{jk})$ are
\[
\alpha_{1k}=\left\{\begin{array}{ll} -\frac{\hat{a}_{k}
b_{k}^1}{\hat{a}_{\mu_1} b_{\mu_1}^1}\,{\cal U}_1^{1} {\cal
U}_{k}^{1} (\omega^1)^{T} & k=1,\ldots,\ell_1\\
\\
-\frac{\hat{a}_{k} b_{k}^1}{\hat{a}_{\mu_1} b_{\mu_1}^1}\,({\cal
U}_1^{1})^2(\omega^1)^{T} & k=1+\ell_1,\ldots,n-1 \\
\\
\frac{\hat{a}_{\mu_1} b_{\mu_1}^1}{\hat{a}_{\mu_1}
b_{\mu_1}^1}({\cal U}_1^{1})^2(\omega^1)^{T} & k=n.
\end{array}
\right.
\]
for $j=2,\ldots,r-1$
\[
\alpha_{jk}=\left\{\begin{array}{ll} -\frac{\hat{a}_{k}
b_{k}^j}{\hat{a}_{\mu_j} b_{\mu_j}^j}\,({\cal
U}_1^{j})^2(\omega^j)^{T} & k=1,\ldots,\mu_{j-1} \AND k=\mu_{j}+1,\ldots,n-1 \\
\\
-\frac{\hat{a}_{k} b_{k}^j}{\hat{a}_{\mu_j} b_{\mu_j}^j}\,{\cal
U}_1^{j} {\cal
U}_{k-\mu_{j-1}}^{j} (\omega^j)^{T} & k=1+\mu_{j-1},\ldots,\mu_{j}\\
\\
\frac{\hat{a}_{\mu_j} b_{\mu_j}^j}{\hat{a}_{\mu_j}
b_{\mu_j}^j}({\cal U}_1^{j})^2(\omega^j)^{T} & k=n.
\end{array}
\right.
\]
and finally
\[
\alpha_{rk}=\left\{\begin{array}{ll} -\frac{\hat{a}_{k}
b_{k}^r}{\hat{a}_{\mu_r} b_{\mu_r}^r}\,{\cal
U}_{\ell_r}^{r} {\cal U}_{1}^{r} (\omega^r)^{T} & k=1,\ldots,\mu_{r-1}\\
\\
-\frac{\hat{a}_{k} b_{k}^r}{\hat{a}_{\mu_r} b_{\mu_r}^r}\,{\cal
U}_{\ell_r}^{r} {\cal
U}_{k-\mu_{r-1}}^{r} (\omega^r)^{T} & k=1+\mu_{r-1},\ldots,\mu_{r}-1\\
\\
\frac{\hat{a}_{\mu_r} b_{\mu_r}^r}{\hat{a}_{\mu_r}
b_{\mu_r}^r}({\cal U}_{\ell_r-1}^{r})^2(\omega^r)^{T} & k=n.
\end{array}
\right.
\]
where we recall that $\mu_{r}=n$. Note that the elements of the last
column are rewritten as to lead to the significant simplification of
the determinant to the following form:
\begin{equation}
{\cal
T}(\hat{\Phi})=\dfrac{(-1)^{n-1}\omega_1^1\cdots\omega_{\ell_1}^1\cdots\omega_1^r\cdots\omega_{\ell_r}^r
\hat{a}_1\cdots\hat{a}_{n-1}}{(\hat{a}_{\mu_1}b_{\mu_1}^1)^{\ell_1}\cdots
(\hat{a}_{\mu_r}b_{\mu_r}^r)^{\ell_r}} \det(\diag({\cal U}_1^1,{\cal
U}_1^2,\ldots,{\cal U}_1^{r-1},{\cal U}_{\ell_r}^r))\det {\cal
I}'_{B}
\label{Tdef-2}
\end{equation}
where
\[
{\cal I}'_{B}=\left(\begin{array}{c} Q_1 \\ Q_2 \\
\vdots\\ Q_r\end{array}\right)
\]
and for $j=1,\ldots,r$
\[
Q_j=\left[
\begin{array}{cccccccccc}
b_{1}^{j} & \cdots & b_{\mu_{j-1}}^{j} & b_{1+\mu_{j-1}}^{j} &  b_{2+\mu_{j-1}}^{j} & \cdots &  b_{\mu_j}^{j} & b_{\mu_{j}+1}^{j} & \cdots & \hat{a}_{\mu_j} b_{\mu_j}^{j}\\
b_{1}^{j} & \cdots & b_{\mu_{j-1}}^{j} & b_{1+\mu_{j-1}}^{j} &  b_{2+\mu_{j-1}}^{j} & \cdots & - b_{\mu_j}^{j} & b_{\mu_{j}+1}^{j} & \cdots & \hat{a}_{\mu_j} b_{\mu_j}^{j}\\
\vdots & \cdots &\vdots & \vdots &\vdots & \cdots & \vdots & \vdots &\cdots & \vdots \\
b_{1}^{j} & \cdots & b_{\mu_{j-1}}^{j}& b_{1+\mu_{j-1}}^{j} &
-b_{2+\mu_{j-1}}^{j} & \cdots & -b_{\mu_j}^{j} & b_{\mu_j+1}^{j} &
\cdots & \hat{a}_{\mu_j} b_{\mu_j}^{j}
\end{array}\right]
\]
is a $\ell_j\times n$ matrix. Moreover, $\det {\cal I}'_{B}\neq 0$
since $\det {\cal I}'_{B}=\hat{a}_n \det{\cal I}_{B}$ and
$\hat{a}_n\neq 0$. Thus ${\cal T}(\hat{\Phi})\neq 0$.

\vspace*{0.25in} It follows that there is a neighborhood
$N'\subseteq N$ in which ${\cal T}\neq 0$.  This is equivalent to
saying that for all $\Phi\in N'$, the image of ${\cal S}_{\Phi}$ by
$G_{\Psi}$ in $\mathrm{V}_{\Psi}$ is transverse to the integral
curves of the vector field $\dot{\Psi}=\omega$.

For each $j=1,\ldots,n-1$, the integral curves of the vector field
$\dot{\Phi}^j=\omega$ are dense in the torus $\mathrm{T}$. Thus, for
any $\varepsilon>0$, there is a $\tau_{j,\varepsilon}>0$ and an
$s_{j,\varepsilon}>0$ such that the integral curve segment
\[
\{\,\Phi^j=\tau_j\omega\,\,\,(\mbox{\rm
mod}\,\,\mathrm{T})\,\,\,|\,\,\,\tau_{j,\varepsilon}-s_{j,\varepsilon}<\tau_j<\tau_{j,\varepsilon}+s_{j,\varepsilon}\,\}
\]
is in the $\varepsilon$-ball centered on $\frac{\pi}{2}v_j$ in
$\mathrm{T}$. For $\varepsilon>0$ small enough, the surface
\[
\{\,\Phi=(\tau_1\omega,\tau_2\omega,\ldots,\tau_{n-1}\omega)\,\,(\mbox{\rm
  mod}\,\,\mathrm{V}_{\Phi})\,\,\,|\,\,\,\tau_{j,\varepsilon}-s_{j,\varepsilon}<\tau_j<\tau_{j,\varepsilon}+s_{j,\varepsilon}\,\}
\]
is contained in $N'$ and coincides with the surface ${\cal
  S}_{\Phi^*}$ for
\[
\Phi^*=(\tau_{1,\varepsilon}\omega,\tau_{2,\varepsilon}\omega,\ldots,\tau_{(n-1),\varepsilon}\omega)\,\,\,\,(\mbox{\rm
  mod}\,\,\mathrm{V}_{\Phi}\,).
\]
  Thus, by our previous
  result,
the $n-1$-dimensional
  surface $G_{\Psi}({\cal S}_{\Phi^*})$ is transverse to the integral
  curves of $\dot{\Psi}=\omega$ in $\mathrm{V}_{\Psi}$.  Since these
  integral curves are dense in $\mathrm{V}_{\Psi}$, there are infinitely many intersections
  with $G_{\Psi}({\cal S}_{\Phi^*})$ near
  the point $\widehat{\Psi}=G_{\Psi}(\widehat{\Phi})$.

Let $\stackrel{\circ}{\Psi}\in G_{\Psi}({\cal S}_{\Phi^{*}})$ be
such an intersection point near $\widehat{\Psi}$.  Then there is a
$\stackrel{\circ}{\tau}_n>0$ such that
\[
\stackrel{\circ}{\Psi}=\stackrel{\circ}{\tau_n}\omega\,\,(\mbox{\rm
  mod}\,\,\mathrm{V}_{\Psi}).
\]
Let $\stackrel{\circ}{\Phi}\in {\cal S}_{\Phi^*}$ be such that
$G_{\Psi}(\stackrel{\circ}{\Phi})=\stackrel{\circ}{\Psi}$. Then
there are
$\stackrel{\circ}{\tau}_1>0,\,\stackrel{\circ}{\tau}_2>0,\,\ldots\,,\stackrel{\circ}{\tau}_{n-1}>0$
such that
\[
\stackrel{\circ}{\Phi}=(\stackrel{\circ}{\tau}_1\omega,\stackrel{\circ}{\tau}_2\omega,\ldots,\stackrel{\circ}{\tau}_{n-1}\omega)\,\,(\mbox{\rm
  mod}\,\,\mathrm{V}_{\Phi}).
\]

It follows from (\ref{ift-1}) that
$F(\stackrel{\circ}{\Phi},\stackrel{\circ}{\Psi},G_A(\stackrel{\circ}{\Phi});\omega)=0$,
and by construction, this corresponds to a solution of
(\ref{chareq_matrix-1}). \qed

\section{Conclusion}
We have shown in this paper that $n$ nonresonant eigenvalues on the
imaginary axis can be realized by a scalar delay-differential
equation with $n$ delays. Moreover, the same is true for any
collection of $n$ imaginary eigenvalues in a neighborhood of an
$n$-tuple of nonresonant imaginary eigenvalues. We have also shown
how these results can be applied to non-scalar delay-differential
equations in the context of symmetric delay-differential equations
where the characteristic equation decomposes according to the
isotypic decomposition. We apply our result to delay-coupled $\D_n$
symmetric cell systems with $n$-odd.

There are several ways of extending the main result of our paper.
One question we did not address in this paper is if $n$ nonresonant
nonzero imaginary eigenvalues constitutes an upper bound for the
realizability by a scalar equation with $n$ delay. The case $n=1$ is
one such example since an easy calculation shows that we can have at
most one imaginary eigenvalue on the imaginary axis. It is likely,
but unknown, if this is also true for general $n$.

One may want to study whether $k$ zero eigenvalues in a single
Jordan block and $\ell$ nonresonant nonzero imaginary eigenvalues
can be realized in a scalar delay-differential equation with
$k+\ell$ delays. This problem may be feasible by modifying the proof
of Theorem~\ref{thm:main} since the nonresonance of the $\ell$
eigenvalues is again present. However, we expect the argument used
in this paper to breakdown for $n$ nonzero imaginary eigenvalues
with resonance. We can also study the same problem as in this paper
but for higher dimensional delay equations. One problem would be to
find out if $n$ nonresonant nonzero imaginary eigenvalues can be
realized by a $m$ dimensional system with $k$ delays. For instance,
it is known that a pair of nonzero imaginary eigenvalues can be
realized by a two-dimensional equation with one delay~\cite{CG82}.
In this case $n=2,m=2,k=1$ and so $n=mk$; is it possible to realize
three nonresonant nonzero eigenvalues or does the relationship
$n=mk$ provide a bound to realizability in general?

Another problem which can be studied is whether a restriction in the
class of delay equation can change the realizability. For instance,
the characteristic equation for a general two-dimensional system
with one-delay $\tau$ is
\[
\lambda^2+a\lambda+b\lambda e^{-\lambda\tau}+c+de^{-\lambda \tau}=0
\]
while for a second-order equation with one delay in the feedback
term we must set $b=0$. We know in this case that two nonresonant
nonzero imaginary eigenvalues can be realized by a second-order
equation with a unique delay in the feedback term~\cite{CBOM95}. In
fact, two imaginary eigenvalues with $1:2$ resonance has been found
in such an equation~\cite{CL98}. The obvious question is to see if
$n$ nonresonant (and resonant) nonzero imaginary eigenvalues can be
realized within the class of $n^{th}$ order scalar equations with
one delay in the feedback term.

An extension of our main result in a direction relevant for studying
bifurcations is to find out whether the $n$ nonresonant nonzero
imaginary eigenvalues can be realized by a scalar delay equation
such that the remaining eigenvalues have negative real parts; that
is, the multiple Hopf point lies at the boundary of the stability
region for the equilibrium solution.

Finally, let us mention the case of linear $T$-periodic equations
with $N+1$ delays:
\begin{equation}\label{linear-periodic}
\dot x=\sum_{j=0}^{N} a_{j}(t)x(t-\tau_{j})
\end{equation}
where each $a_{j}(t)$ is a $T$-periodic $n\times n$ matrix.
Hale~\cite{Hale01} states the following open problem:

\vspace{0.1in} \noindent {\em Is it possible to give a precise upper
bound in terms of $N$ on the number of Floquet multipliers
of~(\ref{linear-periodic}) that can have moduli $1$?} \vspace{0.1in}

\noindent If we restrict equation~(\ref{linear-periodic}) to scalar
equations we can pose a possibly simpler problem which is related to
the main result of our paper:

\vspace{0.1in} \noindent {\em Is it possible to realize $N+1$
complex numbers $e^{\pm i\omega_1},\ldots,e^{\pm i\omega_{N+1}}$
with $\omega_1,\ldots,\omega_{N+1}$ positive and rationally
independent as Floquet multipliers of the scalar
equation~(\ref{linear-periodic})?} \vspace{0.1in}

This problem is automatically solved by Theorem~\ref{thm:main} if a
Floquet theorem can be applied to~(\ref{linear-periodic}); that is,
the Floquet exponents of the Floquet multipliers of
equation~(\ref{linear-periodic}) are eigenvalues of a scalar
equation~(\ref{linear-eq1}) with $N+1$ delays. Such a theorem has
not been proved in general, however it may hold true given some
conditions are imposed on~(\ref{linear-periodic}).

\vspace*{0.25in} \noindent {\Large\bf Acknowledgements}

\vspace*{0.2in} This research is partly supported by the Natural
Sciences and Engineering Research Council of Canada in the form of a
Discovery Grant (PLB,VGL).


\begin{thebibliography}{M}
\bibitem{BC94}
J.~B\'elair and S.A.~Campbell.
\newblock Stability and bifurcations of equilibria in a multiple-delayed
differential equation.
\newblock {\em SIAM J.~Appl.~Math.} {\bf 54}, (1994) 1402--1424.

\bibitem{BBL} A.~Beuter, J.~B\'elair and C.~Labrie. Feedback and delays in neurological diseases : a modeling study
using dynamical systems. {\em Bulletin Math. Biology} {\bf 55},
(1993) 525--541.

\bibitem{BGMT-book} A.~Beuter, L.~Glass, M.~Mackey and M.~Titcombe eds. Nonlinear Dynamics in Physiology and Medicine.
{\em Interdisciplinary Applied Mathematics}, {\bf 25}, Springer,
New-York, 2003.

\bibitem{CBOM95} S.A.~Campbell, J.~B\'elair, T.~Ohira and J.~Milton.
Limit cycles, tori, and complex dynamics in a second-order
differential equation with delayed negative feedback. {\em J. Dynam.
Differential Equations} {\bf 7} (1995), 213--236.

\bibitem{CL98} S.A~Campbell, V.G.~LeBlanc. Resonant Hopf-Hopf
interactions in delay differential equations. {\em J. Dynam.
Differential Equations} {\bf 10} (1998), 327--346.

\bibitem{Campbelletal05} S.A. Campbell, Y. Yuan, S. Bungay. Equivariant Hopf bifurcation
in a ring of identical cells with delayed coupling, {\em
Nonlinearity} {\bf 18} (2005) 2827--2846.

\bibitem{Choi-LeBlanc1} Y.~Choi and V.G. LeBlanc. Toroidal normal forms for
bifurcations in retarded functional differential equations. I.
Multiple Hopf and transcritical/multiple Hopf interaction. {\em J.
Differential Equations} {\bf 227} (2006), 166--203.

\bibitem{CG82} K. Cooke and Z. Grossman. Discrete delay, distributed delay
and stability switches. {\em J.~Math.~Anal.~Appl.} {\bf 86} (1982),
592--627.

\bibitem{Diekmann-etal} O.~Diekmann, S.A.~van Gils, S.M.~Verduyn-Lunel and
H.O.~Walther. Delay-equations, functional, complex, and nonlinear
analysis. {\em Applied Mathematical Sciences} {\bf 110},
Springer-Verlag, New-York, (1995).

\bibitem{Dodlaetal} R. Dodla, A. Sen, G.L. Johnston.
Phase-locked patterns and amplitude death in a ring of delay-coupled
limit cycle oscillators. {\em Phys. Rev. E}. {\bf 69} (2004),
056217.

\bibitem{FMR} T. Faria and L.T. Magalh$\tilde{\mbox{\rm a}}$es. Realisation of ordinary differential equations
by retarded functional differential equations in neighborhoods of
equilibrium points. {\em Proc. Roy. Soc. Ed.} {\bf 125A}, (1995)
759--776.

\bibitem{GM} L.~Glass and M.C.~Mackey. Oscillations and chaos in physiological
control systems. {\em Science} {\bf 197}, (1977) 287--289.

\bibitem{GSS88} M. Golubitsky, I. Stewart and D.G. Schaeffer. Singularities and
Groups in Bifurcation Theory, Vol. II. {\em Applied Mathematical
Sciences}, {\bf 69}, Springer-Verlag, New York, 1988.

\bibitem{Guo05} S. Guo, Spatio-temporal patterns of nonlinear oscillations in an
excitatory ring network with delay. {\em Nonlinearity} {\bf 18}
(2005), 2391--2407.

\bibitem{Guo-Huang1} S. Guo, L. Huang, Hopf bifurcating periodic orbits in a ring of
neurons with delays, Phys. D 183 (2003) 19--44.

\bibitem{Guo-Huang2} S. Guo and L. Huang. Stability of nonlinear waves in a ring of neurons with delays. {\em J. Diff.
Eq} {\bf 236} (2007), 343--374.

\bibitem{Gurney} M.S. Gurney, S.P. Blythe and R.M. Nisbee. Nicholson's blowflies revisited, {\em Nature}
{\bf 287} (1980), 17--21.

\bibitem{Hale01} J.K. Hale. Some problems in FDE. {\em Fields Institute
Communications} {\bf 29} (2001), 195--222.

\bibitem{HL} J.K.~Hale and S.M.~Verduyn-Lunel. Introduction to
functional differential equations. {\em Applied Mathematical
Sciences} {\bf 99}, Springer-Verlag, New-York, (1993).

\bibitem{Kuang} Y.~Kuang.
\newblock {\em Delay differential equations with applications in population
dynamics.}
\newblock Mathematics in Science and Engineering, 191. Academic Press,
Boston, (1993).

\bibitem{Kuznetsov} Y. Kuznetsov. Element of Applied
Bifurcation Theory, 3rd edition. {\em Applied Mathematical Science}
{\bf 112}, Springer, New-York, (2004).

\bibitem{LK} R. Lang and K. Kobayashi. External optical feedback effects on semiconductor injection laser properties.
{\em IEEE J. Quantum Electronics} {\bf 16} (1980), 347--355.

\bibitem{Leite-Golubitsky06} M.~Leite, M.~Golubitsky. Homogeneous three-cell networks. {\em Nonlinearity} {\bf 19} (2006),
2313--2363.

\bibitem{Peng07} M. Peng. Bifurcation and stability analysis of nonlinear waves in $\D_n$-symmetric delay-differential systems.
{\em J. Diff. Eq} {\bf 232} (2007) 521--543.

\bibitem{Peng-Yuan1} M. Peng, Y. Yuan, Complex dynamics in discrete delayed models
with $\D_4$-symmetry. {\em Chaos Solitons Fractals} in press.

\bibitem{RLL02} B.F.~Redmond, V.G.~LeBlanc and A.~Longtin. Bifurcation
analysis of a class of first-order nonlinear delay-differential
equations with reflectional symmetry. {\em Physica D} {\bf 166}
(2002), 131--146.

\bibitem{SA}
E.~Stone and A.~Askari.
\newblock Nonlinear models of chatter in drilling processes.
\newblock {\em Dyn. Syst.} {\bf 17} (2002), 65--85.

\bibitem{SS} M.J. Suarez and P.L. Schopf.
\newblock A Delayed Action Oscillator for ENSO.
\newblock
{\em J. Atmos. Sci.}  {\bf 45} (1988), 3283--3287.


\bibitem{Wu98} J. Wu. Symmetric functional-differential equations and neural
networks with memory. {\em Trans. Amer. Math. Soc.} {\bf 350}
(1998), 4799--4838.

\bibitem{Wuetal} J. Wu, T. Faria, Y.S. Huang. Synchronization and stable
phase-locking in a network of neurons with memory. {\em Math.
Comput. Modelling} {\bf 30} (1999), 117--138.
\end{thebibliography}
\end{document}